\documentclass[12pt]{elsarticle}
\usepackage{amsmath}
\usepackage{graphicx}
\usepackage{wrapfig}
\usepackage{amsmath}
\usepackage{times}
\usepackage{amssymb}
\usepackage{array}
\usepackage{latexsym}
\usepackage{amsthm}
\usepackage{amsfonts}
\usepackage{nomencl}
\usepackage{lineno}

\makenomenclature
\usepackage{hyperref}
\usepackage{array,tabularx}
\usepackage{scalerel,stackengine}
\usepackage{graphics}
\usepackage[english]{babel}
\allowdisplaybreaks

\stackMath
\numberwithin{equation}{section}
\date{}
\newtheorem{t1}{Theorem}[section]

\newtheorem{l1}{Lemma}[section]

\usepackage[utf8]{inputenc}
\allowdisplaybreaks

\begin{document}
\begin{frontmatter} 
\title{Existence and regularity of solutions of a supersonic-sonic patch arising in axisymmetric relativistic transonic flow with general equation of state}

\author{Rahul Barthwal}
\author{T. Raja Sekhar}
\address{Department of Mathematics, Indian Institute of Technology Kharagpur, Kharagpur,  India} 
{
\begin{abstract} 
{In this article, we prove the existence and regularity of a smooth solution for a supersonic-sonic patch arising in a modified Frankl problem in the study of three-dimensional axisymmetric steady isentropic relativistic transonic flows over a symmetric airfoil. We consider a general convex equation of state which makes this problem complicated as well as interesting in the context of the general theory for transonic flows. Such type of patches appear in many transonic flows over an airfoil and flow near the nozzle throat. Here the main difficulty is the coupling of nonhomogeneous terms due to axisymmetry and the sonic degeneracy for the relativistic flow. However, using the well-received characteristic decompositions of angle variables and a partial hodograph transformation we prove the existence and regularity of solution in the partial hodograph plane first. Further, by using an inverse transformation we construct a smooth solution in the physical plane and discuss the uniform regularity of solution up to the associated sonic curve. Finally, we also discuss the uniform regularity of the sonic curve.}
\end{abstract}
}
\end{frontmatter} 
%\linenumbers
{
\section{Introduction}
The transonic flow problems are one of the most important problems in mathematical fluid dynamics since transonic flow appears in various important physical phenomena. In the context of transonic flow problems, the study of supersonic bubbles is of utter importance. For a compressible flow passing the duct, Courant and Friedrichs in their famous book \cite{courant0} described that if the Mach number of the flow is not much below one, then the flow becomes supersonic somewhere on the surface of the duct due to the convexity of the duct and is again purely subsonic throughout the exit section. Similar situations arise naturally in many engineering and aerospace applications, such as the flow over an airfoil or in a flow through an axisymmetric nozzle; see Figure \ref{fig:1a}. We refer readers to the monographs of Bers \cite{bers2016mathematical}, Kuz'min \cite{kuz2003boundary} and Shapiro \cite{shapiro1953dynamics} for more details on transonic flows.

\begin{figure}
    \centering
    \includegraphics[width=4 in]{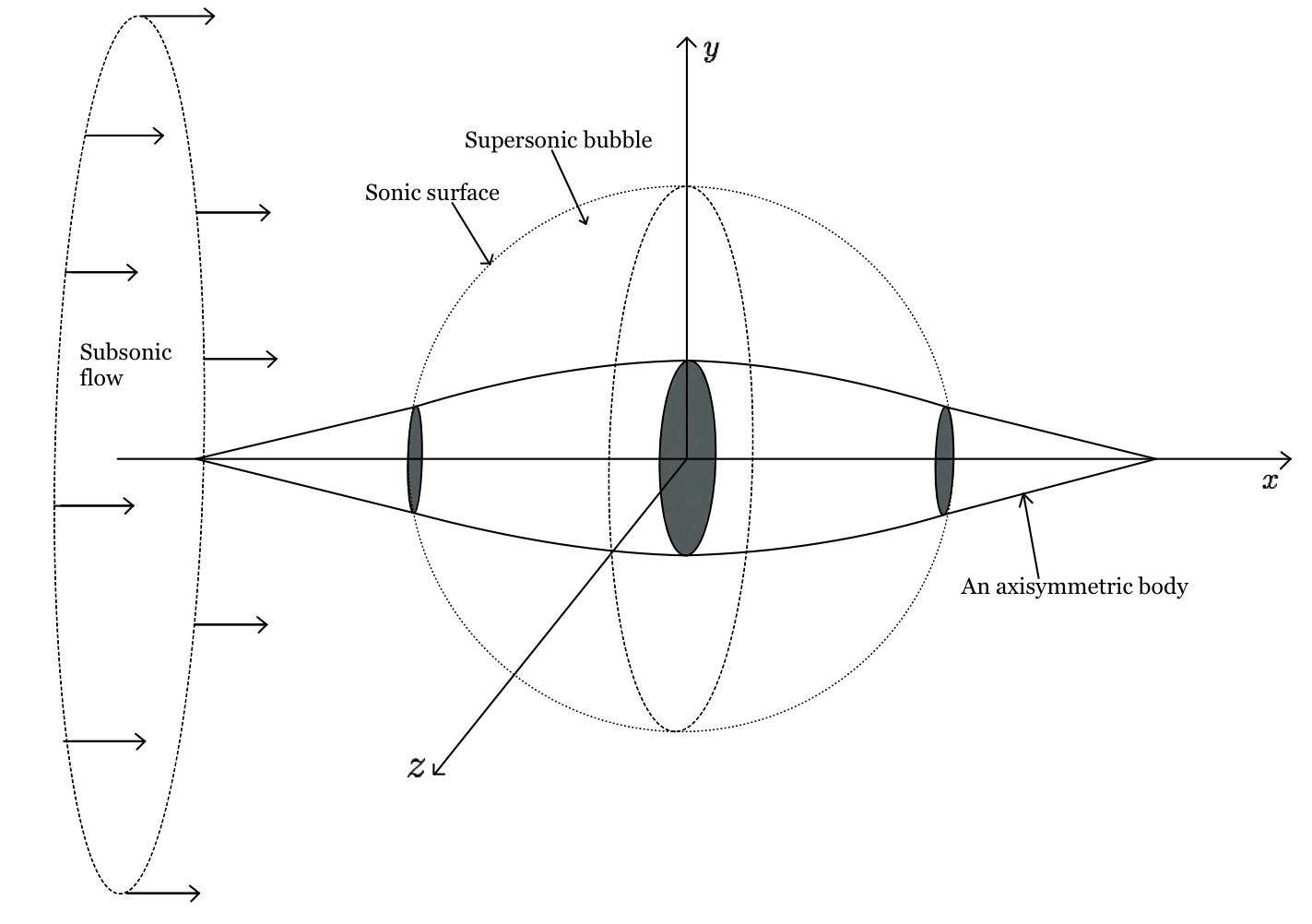}
    \caption{Transonic flow over an axisymmetric body or airfoil}
    \label{fig:1a}
\end{figure}

In the last century, a large number of significant contributions have been made in order to prove the existence of the global transonic solution to such transonic flow problems, but it remains an open mathematical problem till now. The main complexity of the transonic flow is that a transonic structure consists of subsonic and supersonic parts, which are separated by either a sonic curve or transonic shock. These are usually free boundaries due to the nonlinearity of the governing system. Not only this but also the governing systems of transonic flows can change their behavior across the sonic boundary and are usually linearly degenerate on the sonic curve; see \cite{li1998two, li2009interaction, li2011characteristic}. Such features of transonic flow are more complicated to handle when compared to a study of purely subsonic or supersonic flow.

A lot of important existence results for the subsonic-sonic part of the transonic flow for steady Euler equations have been developed in the recent years. Gilbarg and Serrin \cite{gilbarg1955uniqueness} provided a uniqueness result for a subsonic flow past an axisymmetric body, while Xie and Xin \cite{xie2007global, xie2010global} proved the existence of global subsonic-sonic solutions for a 3-D axially symmetric nozzle. In \cite{chen2016subsonic}, Chen et al. established the global existence of a subsonic-sonic solution for the full Euler equations using the compensated-compactness framework. Recently, Wang and Xin \cite{wang2019smooth} proved the existence and uniqueness of a solution for smooth transonic flows of Meyer type in de Laval nozzles and also obtained the first result on the well-posedness for general subsonic-sonic flow problems in \cite{wang2021regular}. On the other hand, for the supersonic part a sonic-supersonic solution for steady isentropic Euler equations was constructed by Zhang and Zheng \cite{zhang2014sonic} while Hu and Li proved the existence of a sonic-supersonic solution for 2-D steady and pseudo-steady full Euler equations \cite{hu2019sonic, hu2020sonic}. The partial hodograph transformation used in the works of Hu and Li; viz. \cite{li2019degenerate, hu2019sonic, hu2020global} has become very crucial while solving sonic-supersonic boundary value problems. We refer readers to \cite{du2011subsonic, du2014steady, chen2016two, wang2013degenerate, hu2021sonic, li2019degenerate, chen2007two} and references cited therein for more such results in the context of sonic-subsonic and sonic-supersonic boundary value problems.
\begin{figure}
    \centering
    \includegraphics[width= 4 in]{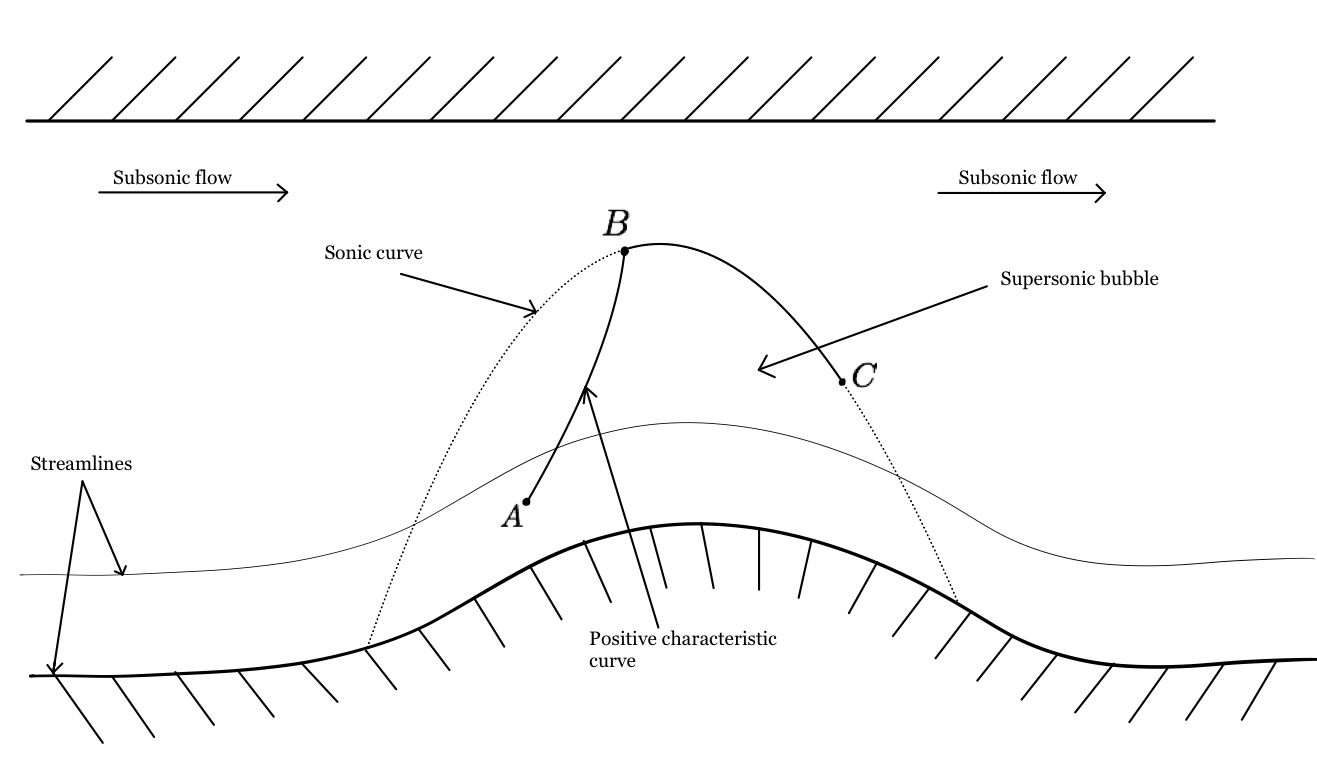}
    \caption{Transonic flow in a channel with a supersonic bubble}
    \label{fig: 1}
\end{figure}

Morawetz, in his work on transonic flow in a channel or a duct, (see Figure \ref{fig: 1}) indicated that a smooth transonic flow does not exist in general, which means that there may exist a transonic shock in the downstream flow \cite{morawetz1964non}. However, it is of utter importance to construct shock-free transonic flows. In \cite{frankl1950formation}, Frankl explored the transonic flow over a symmetric airfoil and suggested that a smooth transonic flow may exist if the part of the airfoil is free of boundary conditions. The original Frankl problem is formulated to find the smooth airfoil's arc $\widehat{EG}$ when the slip conditions on the arcs $\widehat{PE}$ and $\widehat{GQ}$ are prescribed; see Figure \ref{fig: 2}. Many existence and uniqueness results for this Frankl problem have been discussed in the last century; see viz. \cite{morawetz1954uniqueness, cook1978uniqueness}. Kuz'min \cite{kuz2003boundary} proposed a modified Frankl problem in which a velocity distribution is prescribed on the arcs $\widehat{PE}$ and $\widehat{GQ}$ instead of the slip boundary conditions. From a physical point of view, such problems describe the transonic flows past permeable boundaries. The modified Frankl problem can be utilized in many industrial applications as well, where the design of the airfoil and wing usually needs to be formulated according to some specific requirements of the aircraft; for example, one may require the wing profile to have a particular velocity distribution, a specific lift distribution or a certain temperature or pressure distribution. Such a method in the area of aircraft design is usually known as the inverse design method which was pioneered by Lighthill \cite{lighthill1945new} and then developed further by many researchers working in this field; see viz. \cite{hassan1981transonic, henne1981inverse, volpe1981role, volpe1986design, stanitz1988review, labrujere1993computational, obayashi1996genetic} and references cited therein.

\begin{figure}
    \centering
    \includegraphics[width= 4.2 in]{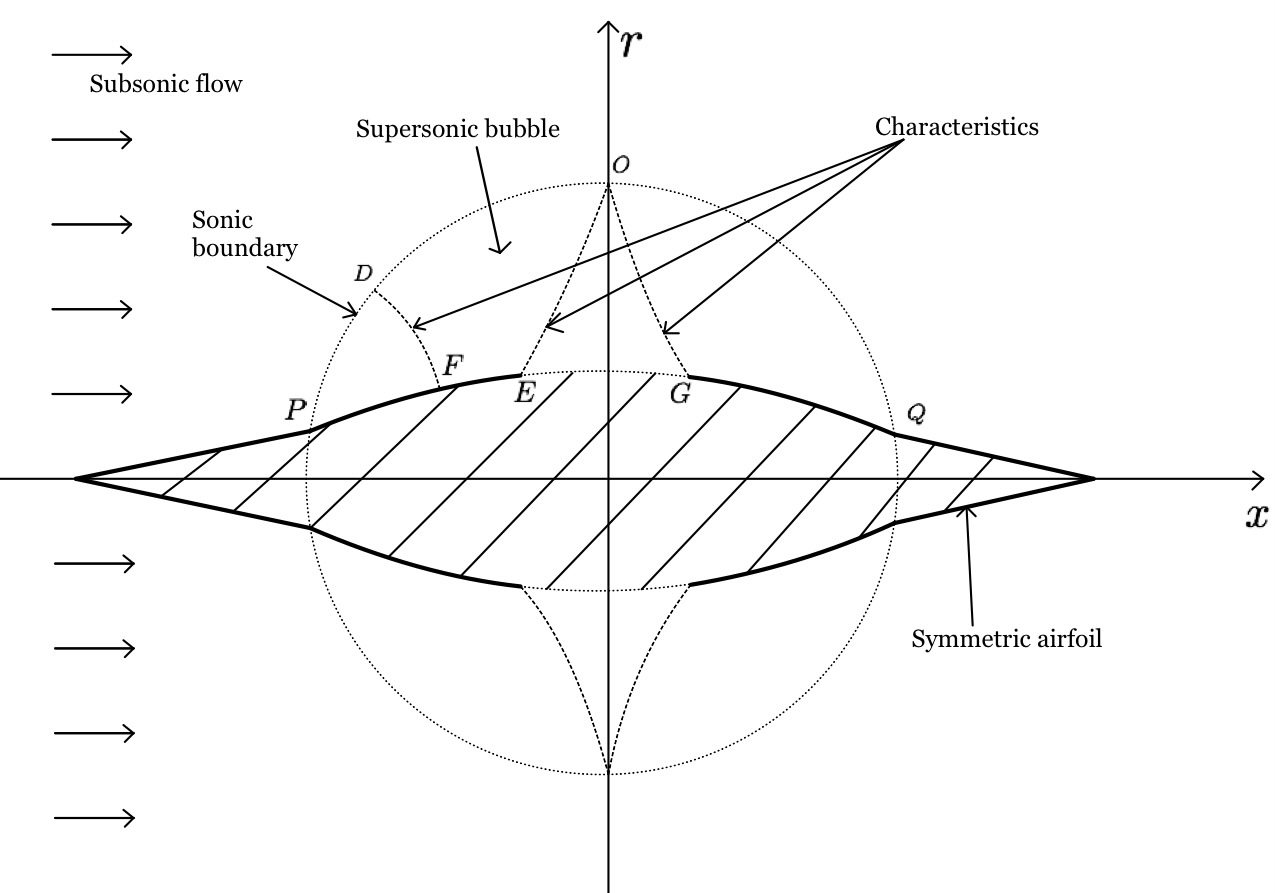}
    \caption{Modified Frankl problem for a transonic flow over a symmetric airfoil: If the velocity distribution on the arcs $\widehat{PE}$ and $\widehat{GQ}$ are prescribed, then find the airfoil's arc $\widehat{EG}$ free of boundary conditions for the correctness of the problem in the class of smooth solutions}
    \label{fig: 2}
\end{figure}

The modified Frankl problem has been studied extensively in the recent past. Kuz'min discussed the existence and uniqueness of the solution of the modified Frankl problem for a linearized version of the von Karman equation in a finite domain; see \cite{kuz2001solvability, kuz2004modified}. Recently, Hu and Li established the existence and regularity of solutions of a sonic-supersonic patch extracted from a modified Frankl problem for 2-D steady isentropic Euler equations and 3-D steady axisymmetric isentropic Euler equations with ideal gas \cite{husonic2021, hu2022supersonic}. The recent development in the context of transonic flows has motivated us to ask naturally whether such analysis can be performed for more complicated mixed-type systems for a more general equation of state or not. Inspired by this idea, the main motivation to do this work is to develop the existence and regularity of a sonic-supersonic solution arising in a modified Frankl problem for 3-D steady axisymmetric relativistic Euler equations with a general convex equation of state. 

In the domain of astrophysics, plasma physics, and nuclear physics, the velocity of fluid particles are usually very large and often very close to the speed of light as well, which means that the relativistic effects have to be taken into consideration and the classical Euler equations of gas dynamics are no longer valid. In the case of such high-speed flow, the governing system under consideration is referred to as the relativistic Euler system. In the recent years, a lot of interesting work has been done in the context of relativistic gas dynamics systems; see viz. \cite{luan2018two, li2005global, chen2018boundary} and references cited therein. Here, we consider the three-dimensional steady isentropic relativistic Euler equations of the form
\begin{equation}\label{eq: 1.a}
\begin{cases}
\begin{aligned}
     &(n\gamma u_1)_x+(n\gamma u_2)_y+(n\gamma u_3)_z=0,\\
     &((\rho+p)\gamma^2 u_1^2+p)_x+((\rho+p)\gamma^2 u_1 u_2)_y+((\rho+p)\gamma^2 u_1 u_3)_z=0,\\
     &((\rho+p)\gamma^2 u_1u_2)_x+((\rho+p)\gamma^2 u_2^2+p)_y+((\rho+p)\gamma^2 u_2u_3)_z=0,\\
     & ((\rho+p)\gamma^2 u_1u_3)_x+((\rho+p)\gamma^2 u_2 u_3)_y+((\rho+p)\gamma^2 u_3^2+p)_z=0,
\end{aligned}
\end{cases}
\end{equation}
where $n$ denotes the proper number density, $u_1, u_2$ and $u_3$ are the velocity components of the velocity along $x, y$ and $z$ directions, respectively. $\gamma=\left(\sqrt{1-u_1^2-u_2^2-u_3^2}\right)^{-1}~(0<u_1^2+u_2^2+u_3^2<1)$ is the Lorentz factor such that the speed of light is normalized to be 1, $\rho+p(\rho)=i$ denotes the enthalpy per unit volume and $p=p(\rho)$ is the pressure of the relativistic fluid with $\rho$ being the total mass-energy density. In the cylindrical coordinates $(x, r, \sigma)$ a flow is said to be axisymmetric if the state variables are independent of the angle $\sigma$. Further, if we consider that the flow is axisymmetric about the $x$ axis and is without swirl, i.e., $(\rho, u_1, u_2, u_3)$ satisfy
\begin{align*}
\begin{cases}
\rho(x, y, z)=\rho(x, r),~~n(x, y, z)=n(x, r),~~u_1(x, y, z)=u(x, r),\\
u_2(x, y, z)=v(x, r)\cos \sigma, ~~u_3(x, y, z)=v(x, r)\sin \sigma,
\end{cases}
\end{align*}
where $u$ and $v$ are the axial and radial velocity components, respectively. 

The system \eqref{eq: 1.a} can be now rewritten in terms of $(\rho, u, v)(x, r)$ as follows: 
\begin{align}\label{eq: 1.1}
\begin{cases}
    &(n\gamma u)_x+(n\gamma v)_r=-\dfrac{n\gamma v}{r},\\
     &((\rho+p)\gamma^2 u^2+p)_x+((\rho+p)\gamma^2 uv)_r=-\dfrac{(\rho+p)\gamma^2 uv}{r},\\
     &((\rho+p)\gamma^2 uv)_x+((\rho+p)\gamma^2 v^2+p)_r=-\dfrac{(\rho+p)\gamma^2 v^2}{r},
\end{cases}     
\end{align}
where $\gamma=(1-q^2)^{-1/2}$ is the normalized Lorentz factor of axisymmetric relativistic flow and $0<q=\sqrt{u^2+v^2}<1$ is the flow velocity. Further, throughout the article we assume that the mass-energy density $\rho$ and pressure $p=p(\rho)$ satisfies \cite{chen2018boundary, chen2004stability}
\begin{align}\label{1.2}
    0<\rho<\rho_{\max}<\infty,~~ 0<p'(\rho)<1, ~~p''(\rho)>0.
\end{align}
and $p''(\rho)$ remains finite for all values of $\rho$, which is generally true for all physically relevant equations of states.

Now noting the Figure \ref{fig: 2}, we define the sonic-supersonic problem under consideration for 3-D axisymmetric relativistic Euler equations precisely as follows:\\\\
\textbf{Supersonic-sonic boundary value problem extracted from modified Frankl problem for 3-D relativistic flow:\\}
\textit{If $\widehat{PE}$ is an increasing and concave smooth streamline of an axisymmetric relativistic transonic flow and a velocity distribution is prescribed on the arc $\widehat{PE}$ such that the point $P$ is sonic, then find a sonic curve $\widehat{PD}$ starting from point $P$ and construct a smooth supersonic solution for 3-D axisymmetric steady isentropic relativistic Euler equations in a region $PFD$ near the point $P$ bounded by the sonic curve $\widehat{PD}$, the arc $\widehat{PF}$ and a negative characteristic curve $\widehat{DF}$ for a general convex pressure. Moreover, check the regularity of the constructed solution.}\\

One of the main complexity of the problem under consideration is that the velocity data is given only on the streamline arc $\widehat{PE}$ in contrast to all other sonic-supersonic boundary value problems or semi-hyperbolic patch problems where the data is prescribed not only on a streamline but also on a characteristic curve as well (see for example \cite{hu2020global, li2019degenerate, li2011semi, barthwal2022existence}). In particular, for relativistic Euler equations, we refer readers to \cite{fan2022sonic}. The other important complexity in this problem is to handle the nonhomogeneous terms due to the axisymmetry and the sonic degeneracy along the sonic curve. In all the previous work related to 2-D steady systems, angle variables (Mach angle and flow angle) were chosen as independent variables to convert the governing system into a linearized one. However, one can not expect to linearize the axisymmetric systems due to the presence of nonhomogeneous terms. To overcome these complexities, we use partial hodograph transformation where the independent variables are Mach angle and the potential function to convert the governing axisymmetric relativistic Euler system into a new degenerate hyperbolic system. The idea of choosing such independent variables is taken from a very recent work of Hu \cite{hu2022supersonic}. However, unlike the Mach-flow angle plane, the reduced hyperbolic equations in our case do not form a closed system and additional equations are needed to be added to the system in order to close the system which makes the current problem even more complicated. We also comment that the derivation of a priori estimates of solutions for the current problem is also not very easy as a priori estimates developed in the previous works such as \cite{lai2015centered, barthwal2022existence, li2009interaction, li2011characteristic, sheng2018interaction} which are based on characteristic decompositions in homogeneous form. But the nonhomogeneous terms in this problem lead us to the nonhomogeneous form of characteristic decompositions of the angle variables, which greatly affect the establishment of a priori estimates of the solutions. However, using some proper auxiliary functions and characteristic decompositions on them, we are able to develop the $C^0$ and $C^1$ estimates of the solutions of this new degenerate hyperbolic system in the partial hodograph plane, which helps us to develop a global solution and its regularity in the partial hodograph plane. Finally, using an inverse transformation, we transform these solutions back to the physical plane in order to solve the original problem. 

The rest of the article can be organized in the following manner. In section 2, we discuss the basic properties of the axisymmetric steady isentropic relativistic Euler equations \eqref{eq: 1.1} and define the characteristic angles for relativistic flow. Section 3 is devoted to defining the problem precisely and prescribing the boundary data on the arc $\widehat{PE}$. Using a partial hodograph transformation, we discuss the existence and regularity of solutions to sonic-supersonic boundary value problem in the new coordinate system in section 4. In section 5, we transform the constructed solutions back into the physical plane by using an inverse transformation and verify that the solutions constructed actually satisfy the boundary value problem. Finally, we provide conclusions and the future scope of this work in section 6.
}
\section{Basic properties of three-dimensional axisymmetric isentropic irrotational steady relativistic Euler equations}\label{2}
We assume that the relativistic flow is irrotational, i.e., $u_r=v_x$ then by first equation of \eqref{eq: 1.1} and making use of $\gamma_x=\gamma^3qq_x, ~~\gamma_r=\gamma^3qq_r,$ second equation of system \eqref{eq: 1.1} becomes
\begin{equation}\label{eq: 2.3}
    \dfrac{n}{\gamma}\bigg\{\dfrac{i}{n}\gamma_x+\dfrac{\gamma}{n}p_x+\dfrac{\gamma u}{n}\bigg[n\gamma u\left(\dfrac{\gamma i}{n}\right)_x+n\gamma v\left(\dfrac{\gamma i}{n}\right)_r\bigg]\bigg\}=0.
\end{equation}
Further, by the second law of thermodynamics, we have
\begin{align}\label{eq: 2.4}
    d\left(\dfrac{i}{n}\right)=\dfrac{1}{n}dp+Tds,
\end{align}
where $T$ is the absolute temperature and $s$ is the entropy of the flow. Since the flow is assumed to be isentropic, i.e.,  $s$ is constant or in other words $ds=0$. Therefore, in view of \eqref{eq: 2.3} and \eqref{eq: 2.4}, we have
\begin{align}\label{eq: a}
    \left(\dfrac{\gamma i}{n}\right)_x+\dfrac{\gamma u}{n}\bigg[n\gamma u\left(\dfrac{\gamma i}{n}\right)_x+n\gamma v\left(\dfrac{\gamma i}{n}\right)_r\bigg]=0.
\end{align}
Similarly, from the third equation of system \eqref{eq: 1.1}, one can easily obtain 
\begin{align}\label{eq: b}
    \left(\dfrac{\gamma i}{n}\right)_r+\dfrac{\gamma v}{n}\bigg[n\gamma u\left(\dfrac{\gamma i}{n}\right)_x+n\gamma v\left(\dfrac{\gamma i}{n}\right)_r\bigg]=0.
\end{align}
It is easy to see that \eqref{eq: a} and \eqref{eq: b} form a homogeneous system of linear equations for $\left(\dfrac{\gamma i}{n}\right)_x$ and $\left(\dfrac{\gamma i}{n}\right)_r$. Now the determinant of the coefficient matrix is 
\begin{align*}
    \begin{vmatrix}
    1+\gamma^2u^2&\gamma^2 uv\\
    \gamma^2uv&1+\gamma^2v^2
    \end{vmatrix}=\dfrac{1}{1-q^2}\neq 0.
\end{align*}
Hence, we must have $\left(\dfrac{\gamma i}{n}\right)_x=\left(\dfrac{\gamma i}{n}\right)_r=0$, which provides the Bernoulli's law for axisymmetric steady relativistic Euler equations of the form
\begin{align}\label{eq: 2.7}
\dfrac{\gamma i}{n}= \mathrm{const.}
\end{align}
For the convenience of the subsequent discussion, we write Bernoulli's law in the following form
\begin{align}\label{eq: 2.8}
\dfrac{\gamma i}{n}= m\hat{\gamma},
\end{align}
where $m$ is the average rest mass per particle and $\hat{\gamma}^{-1}=\sqrt{1-\hat{q}^2}$ is a constant.
\begin{l1}
If $p$ satisfies $\dfrac{\partial p}{\partial n}>0$ for $n>0$ then there exists a constant $\hat{q}~(0<\hat{q}<1)$ such that the flow speed $q<\hat{q}$. The quantity $\hat{q}$ is called the limit speed of the flow and the flow speed approaches the limit speed when $n$ approaches $0$ \cite{chen2018boundary}.
\end{l1}
\begin{proof}
Using the second law of thermodynamics for isentropic flow, we have
$$\dfrac{d\left(\frac{i}{n}\right)}{dn}=\dfrac{d\left(\frac{i}{n}\right)}{dp}.\dfrac{dp}{dn}=\dfrac{1}{n}.\dfrac{dp}{dn}>0$$ 
for $n>0$.

Therefore, using the Bernoulli's law $\eqref{eq: 2.8}$ and the fact that $\dfrac{i}{n}=\dfrac{mn+e+p}{n}\geq m$ for $n\geq 0$, it is easy to see that $q<\hat{q}<1$ and $q$ approaches $\hat{q}$ as $n$ approaches $0$.
\end{proof}
Now noting the Bernoulli's law \eqref{eq: 2.8}, system \eqref{eq: 1.1}
can be rewritten as
\begin{align}\label{eq: 2.10 a}
\begin{cases}
    \gamma[(i\gamma u)_x+(i\gamma v)_r]+i\gamma [u\gamma_x+v\gamma_r]=-\dfrac{i\gamma^2 v}{r},\\
    \gamma u[(i\gamma u)_x+(i\gamma v)_r]+i\gamma [u(\gamma u)_x+v(\gamma u)_r]+p_x=-\dfrac{i\gamma^2 uv}{r},\\
    \gamma v[(i\gamma u)_x+(i\gamma v)_r]+i\gamma [u(\gamma v)_x+v(\gamma v)_r]+p_r=-\dfrac{i\gamma^2 v^2}{r},
\end{cases}    
\end{align}
Then by taking the scalar product of \eqref{eq: 2.10 a} with $\gamma(1, -u, -v)$ and simplifying, we obtain
\begin{align}\label{eq: 2.11 a}
    i(\gamma u)_x+i (\gamma v)_r=-\gamma u \rho_x-\gamma v \rho_r-\dfrac{i \gamma v}{r}.
\end{align}
Again from the momentum equations of \eqref{eq: 1.1}, we can easily obtain 
\begin{align}\label{eq: 2.12 a}
    \begin{cases}
        i\gamma^2 u u_x+i \gamma^2 vu_r+p_x=0,\\
        i\gamma^2 u v_x+i \gamma^2 vv_r+p_r=0.
    \end{cases}
\end{align}
Then taking the scalar product of \eqref{eq: 2.12 a} with $\gamma(u, v)$, we have
\begin{align}\label{eq: 2.13 a}
    i\gamma^3[u^2 u_x+uv(u_r+v_x)+v^2v_r]=-a^2(\gamma u \rho_x+\gamma v \rho_r),
\end{align}
where $a=\sqrt{p'(\rho)}$ denotes the speed of sound relative to the fluid.

Then by combining \eqref{eq: 2.11 a} and \eqref{eq: 2.13 a} the three-dimensional axisymmetric steady isentropic irrotational relativistic flow can be governed by Bernoulli's law \eqref{eq: 2.8} and
\begin{align}\label{eq: 2.9}
    \begin{cases}
    (M_1^2-1)u_x+M_1M_2(u_r+v_x)+(M_2^2-1)v_r=\dfrac{v}{r},\\
    u_r-v_x=0,
    \end{cases}
\end{align}
where $M_1=\dfrac{\gamma u}{a\gamma_a}, ~~M_2=\dfrac{\gamma v}{a\gamma_a},$ and $\gamma_a=\dfrac{1}{\sqrt{1-a^2}}$.
 
In matrix form \eqref{eq: 2.9} can be rewritten as
\begin{align}\label{eq: 2.11}
    \begin{bmatrix}
    M_1^2-1&M_1M_2\\
    0&-1
    \end{bmatrix}
    \begin{bmatrix}
    u\\
    v
    \end{bmatrix}_x
    +
    \begin{bmatrix}
    M_1M_2&M_2^2-1\\
    1&0
    \end{bmatrix}
    \begin{bmatrix}
    u\\
    v
    \end{bmatrix}_r=
    \begin{bmatrix}
    \dfrac{v}{r}\\0
    \end{bmatrix}.
\end{align}
 \begin{figure}
    \centering
    \includegraphics[width=4.2 in]{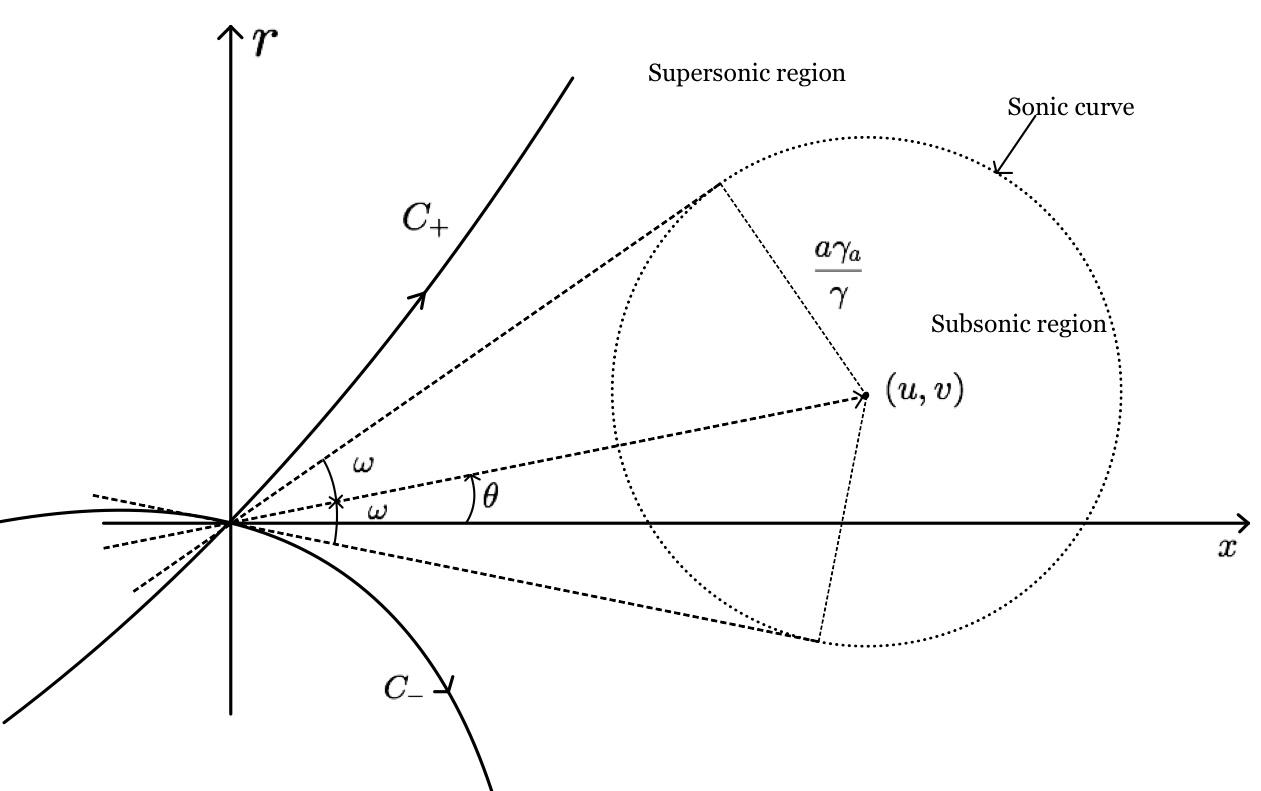}
    \caption{$C_\pm$ characteristic directions, proper Mach angle and flow angle}
    \label{fig: 3}
\end{figure} 

It is easy to see that the system \eqref{eq: 2.11} has the eigenvalues $\Lambda_\pm=\dfrac{M_1M_2\pm \sqrt{M_1^2+M_2^2-1}}{M_1^2-1}$ with corresponding left eigenvectors $l_\pm=(1, \mp \sqrt{M^2-1})$, where $M=\sqrt{M_1^2+M_2^2}=\dfrac{\gamma q}{\gamma_a a}$ is the proper Mach number of the relativistic flow. The expression of these eigenvalues shows that the system \eqref{eq: 2.11} is a mixed-type system and can change its behavior from hyperbolic to elliptic across the sonic boundary ($M=1$); therefore, its behavior depends on the choice of proper Mach number. For $M> 1$ (supersonic) system \eqref{eq: 2.9} is hyperbolic while for $M<1$ (subsonic) it is elliptic. Then we define the two families of wave characteristics as 
\begin{equation}
    \dfrac{dr}{dx}=\Lambda_\pm.
\end{equation}
Moreover, we obtain the characteristic equations by multiplying $l_{\pm}$ to the system $\eqref{eq: 2.11}$ as 
\begin{align}\label{eq: 2.6}
\begin{cases}
{\partial}_{+}u+\Lambda_- {\partial}_{+}v=\dfrac{v}{r(M_1^2-1)},\\
{\partial}_{-}u+\Lambda_+ {\partial}_{-}v=\dfrac{v}{r(M_1^2-1)}
\end{cases}
\end{align} 
where ${\partial}_{\pm}=\partial_x+\Lambda_{\pm}\partial_r$.

From the expression of eigenvalues $\Lambda_\pm$, it is easy to see that
\begin{align}
    1=\dfrac{|(M_1, M_2).(\Lambda, -1)|}{|\Lambda, -1|},
\end{align}
which means that the component of the flow velocity normal to $C_\pm$ characteristic curve is equal to $\dfrac{a\gamma_a}{\gamma}$. 
Then we define the concept of characteristic direction as in \cite{lai2014characteristic}. The direction of the characteristic is defined as the tangential direction that forms an acute angle $\omega$ with the flow velocity vector $(u, v)$. Geometrically, the $C_+$ characteristic direction forms the angle $\omega$ with the flow velocity vector $(u, v)$ in a clockwise direction, while the $C_-$ characteristic direction forms the angle $\omega$ with the flow velocity vector $(u, v)$ in the counterclockwise direction where $\omega$ is called the proper Mach angle. Further, we denote the flow angle by $\theta$, which is the angle between velocity vector $(u, v)$ and $x$-axis such that
$\tan \theta=\dfrac{v}{u}$ and $\sin \omega=\dfrac{1}{M}$
; see Figure \ref{fig: 3}. 
\subsection{Characteristic equations in terms of characteristic angles}
First we differentiate Bernoulli's law \eqref{eq: 2.7} w.r.t. $q$, which yields
\begin{align}\label{eq: 2.14}
    \dfrac{\gamma^3 q i}{n}+\gamma\dfrac{d(i/n)}{dp}\dfrac{dp}{da}\dfrac{da}{dq}=0.
\end{align}
Then from $a^2=p'(\rho)$ we have $\dfrac{dp}{da}=\dfrac{2a^3}{p''(\rho)}$, which is exploited in \eqref{eq: 2.14} along with second law of thermodynamics \eqref{eq: 2.4} to yield
\begin{align}
    \dfrac{da}{dq}=-\dfrac{iq\gamma^2 p''(\rho)}{2a^3}<0.
\end{align}
Also, by $M=\dfrac{\gamma q}{a\gamma_a}$, it is easy to see that
\begin{align}
    \dfrac{dM}{dq}=M\left(\dfrac{1}{q(1-q^2)}-\dfrac{1}{a(1-a^2)}\dfrac{da}{dq}\right).
\end{align}
Then noting that $\dfrac{da}{dq}<0$ and $0<a^2, q^2<1$ we have $\dfrac{dM}{dq}>0$.

 We invoke polar coordinates in velocity plane such that $u=q\cos \theta$ and $v=q\sin \theta$. Then we can use the following formulas of velocity \cite{chen2018boundary} 
 \begin{equation}
 u=\dfrac{a\gamma_a\cos\theta}{\gamma \sin \omega}, ~~v=\dfrac{a\gamma_a\sin\theta}{\gamma \sin \omega}.
 \end{equation}
We further introduce the weighted directional derivatives along the characteristics \cite{hu2021degenerate, hu2022supersonic}
\begin{align}
\begin{cases}
    &\Tilde{\partial}_+=r \cos \alpha \partial_x+r \sin \alpha \partial_r, ~~\Tilde{\partial}_-=r \cos \beta \partial_x+r \sin \beta \partial_r,\\
    &\Tilde{\partial}_0=r \cos \theta \partial_x+r \sin \theta \partial_r,
    \end{cases}
\end{align}
from which one has
\begin{align}
    \partial_x=-\dfrac{\sin \beta \Tilde{\partial}_+-\sin \alpha \Tilde{\partial}_-}{r \sin (2\omega)},~~\partial_r=\dfrac{\cos \beta \Tilde{\partial}_+-\cos \alpha \Tilde{\partial}_-}{r \sin (2\omega)},~~\Tilde{\partial}_0=\dfrac{\Tilde{\partial}_++\Tilde{\partial}_-}{2\cos \omega}.
\end{align}
Then the characteristic equations \eqref{eq: 2.6} can be rewritten as
\begin{align}\label{eq: 2.25a}
    \begin{cases}
        \Tilde{\partial}_+u+\Lambda_- \Tilde{\partial}_+ v=\dfrac{v \cos \alpha}{M_1^2-1},\\
        \Tilde{\partial}_-u+\Lambda_+ \Tilde{\partial}_- v=\dfrac{v \cos \beta}{M_1^2-1},
    \end{cases}
\end{align}
Also, in terms of weighted directional derivatives, we have the first-order decompositions of velocity components as
 \begin{align}\label{eq: 2.5}
     \begin{cases}
     \Tilde{\partial}_\pm u=\dfrac{\gamma_a}{\gamma \sin^2 \omega}\bigg[\sin \omega \cos \theta f(a)\Tilde{\partial}_\pm a-a  \sin \omega \sin \theta \Tilde{\partial}_\pm \theta-a \cos \omega \cos \theta \Tilde{\partial}_\pm\omega\bigg],\vspace{0.2 cm}\\
    \Tilde{\partial}_\pm v=\dfrac{\gamma_a}{\gamma \sin^2 \omega}\bigg[\sin \omega \sin \theta f(a)\Tilde{\partial}_\pm a+a  \sin \omega \cos \theta \Tilde{\partial}_\pm \theta-a \cos \omega \sin \theta \Tilde{\partial}_\pm\omega\bigg],
     \end{cases}
 \end{align}
 where $0<f(a)=\left(\gamma_a^2+\dfrac{2a^4}{ip''(\rho)}\right)<\infty$.
 
 Now using \eqref{eq: 2.5} in characteristic equation \eqref{eq: 2.25a}, we have
 \begin{align}
 \begin{cases}
     \Tilde{\partial}_+ \omega&=\dfrac{aF_1(a, \omega)}{2iq\gamma p''(\rho)\cos \omega \gamma_a}\dfrac{\Tilde{\partial}_+a}{a},\\
     \Tilde{\partial}_- \omega&=\dfrac{aF_1(a, \omega)}{2iq\gamma p''(\rho)\cos \omega \gamma_a}\dfrac{\Tilde{\partial}_-a}{a},\\
     \Tilde{\partial}_+ \theta&=\dfrac{(F_2(a, \omega)-aF_1(a, \omega)\cos^2\omega)}{2iq\gamma p''(\rho)\cos \omega\sin^2\omega \gamma_a}\dfrac{\Tilde{\partial}_+a}{a}-\sin \omega \sin \theta,\\
      \Tilde{\partial}_- \theta&=\dfrac{(aF_1(a, \omega)\cos^2\omega-F_2(a, \omega))}{2iq\gamma p''(\rho)\cos \omega\sin^2\omega \gamma_a}\dfrac{\Tilde{\partial}_-a}{a}+\sin \omega \sin \theta,
      \end{cases}
 \end{align}
 where 
 \begin{align}
 \begin{cases}
 F_1(a, \omega)=2ip''(\rho) \gamma_a^2f(a)+4 a^2\sin^2\omega> 0,\\
 F_2(a, \omega)=iq\gamma\gamma_a p''(\rho)f(a)\cos \omega \sin 2\omega.
 \end{cases}
 \end{align}
 If $\omega=\omega(a)$ then one can use $\dfrac{d\omega}{da}=\dfrac{d\omega}{dM}\dfrac{dM}{dq}\dfrac{dq}{da}$ and $\dfrac{da}{dq}<0,~\dfrac{dM}{dq}>0$ together with $\sin \omega=\dfrac{1}{M}$ to yield $\dfrac{d\omega}{da}>0$. Let $\varpi=\sin \omega, \omega\in [k_0, \pi/2], ~k_0>0$ is a constant. Then by inverse function theorem we must have $a=a(\omega)=a(\sin^{-1}\varpi)$. Therefore we can define $F_1(a, \omega)=:F_1(\varpi)$ and $F_2(a, \omega)=:F_2(\varpi)$ to obtain 
 \begin{align}\label{eq: 2.21}
     \begin{cases}
     \Tilde{\partial}_+\theta+\dfrac{4a^2\cos \omega}{F_1(\varpi)}\Tilde{\partial}_+\varpi=-\varpi \sin \theta,\\
     \Tilde{\partial}_-\theta-\dfrac{4a^2 \cos \omega}{F_1(\varpi)}\Tilde{\partial}_-\varpi=\varpi \sin \theta.
     \end{cases}
 \end{align}
 Due to the continuity of the function $\dfrac{2a^2(y)}{yF_1(y)}, y\in [\sin k_0, 1]$, we set
 \begin{align*}
     I(\varpi)=\int_{\sin k_0}^{\varpi}\dfrac{2a^2(y)}{yF_1(y)}dy, ~~\varpi\in [\sin k_0, 1].
 \end{align*}
 Therefore, we write $I:=I(\varpi)$ to convert \eqref{eq: 2.21} in the following form
 \begin{align}\label{eq: 2.23}
     \begin{cases}
     \Tilde{\partial}_+\theta+\sin 2\omega \Tilde{\partial}_+I=-\varpi \sin \theta,\\
     \Tilde{\partial}_-\theta-\sin 2\omega \Tilde{\partial}_-I=\varpi \sin \theta
     \end{cases}
 \end{align}
 with 
 \begin{align}\label{eq: 2.24}
     \Tilde{\partial}_i \varpi=\dfrac{\varpi F_1(\varpi)}{2a^2(\varpi)}\Tilde{\partial}_i I, ~~i=0, \pm.
 \end{align}
 Now we use the following commutator relation from  \cite{hu2022supersonic}
 \begin{align*}
  \Tilde{\partial}_-\Tilde{\partial}_+ -\Tilde{\partial}_+\Tilde{\partial}_-&=\dfrac{1}{\sin 2\omega}\bigg[(\cos 2\omega \Tilde{\partial}_+\beta-\Tilde{\partial}_-\alpha)\Tilde{\partial}_- -(\Tilde{\partial}_+\beta-\cos 2 \omega \Tilde{\partial}_-\alpha)\Tilde{\partial}_+\bigg]+\sin \beta \Tilde{\partial}^+-\sin \alpha \Tilde{\partial}^-
 \end{align*}
to obtain the commutator relation of $I$ of the form
\begin{align*}
   \Tilde{\partial}_-\Tilde{\partial}_+I-\Tilde{\partial}_+\Tilde{\partial}_-I&=\dfrac{1}{\sin 2\omega}\bigg[(\cos 2\omega \Tilde{\partial}_+\beta-\Tilde{\partial}_-\alpha)\Tilde{\partial}_-I-(\Tilde{\partial}_+\beta-\cos 2 \omega \Tilde{\partial}_-\alpha)\Tilde{\partial}_+I\bigg]+\sin \beta \Tilde{\partial}^+I-\sin \alpha \Tilde{\partial}^-I, \\
   \Tilde{\partial}_-\Tilde{\partial}_+\theta-\Tilde{\partial}_+\Tilde{\partial}_-\theta&=\dfrac{1}{\sin 2\omega}\bigg[(\cos 2\omega \Tilde{\partial}_+\beta-\Tilde{\partial}_-\alpha)\Tilde{\partial}_-\theta-(\Tilde{\partial}_+\beta-\cos 2 \omega \Tilde{\partial}_-\alpha)\Tilde{\partial}_+\theta\bigg]+\sin \beta \Tilde{\partial}^+I-\sin \alpha \Tilde{\partial}^-I.
 \end{align*}
Therefore, if we denote $W=\Tilde{\partial}_+I$, $Z=-\Tilde{\partial}_-I$,  $f(a):=f(a(\varpi))$ and use $F_1(\varpi)=4a^2\varpi^2+2ip''(\rho)\gamma_{a(\varpi)}^2f(a(\varpi))$, then it is easy to obtain the characteristic decompositions of $W$ and $Z$ of the form
 \begin{align}\label{eq: 2.28}
     \begin{cases}
         &\hspace{-0.5 cm}\Tilde{\partial}_-W=W\Bigg[\dfrac{F_1(\varpi)}{4a^2\cos^2 \varpi}(W-Z)+W-\cos 2\omega~ Z+\dfrac{F_1(\varpi)}{2a^2}Z+2\sin \theta \cos \omega\Bigg]+\dfrac{\sin \beta}{2}(W-Z)+\sin^2\theta,\vspace{0.2 cm}\\
    &\hspace{-0.5 cm} \Tilde{\partial}_+Z=Z\Bigg[\dfrac{F_1(\varpi)}{4a^2\cos^2 \varpi}(W-Z)-Z+\cos 2\omega~ W-\dfrac{F_1(\varpi)}{2a^2}W+2\sin \theta \cos \omega\Bigg]-\dfrac{\sin \alpha}{2}(W-Z)-\sin^2\theta,
     \end{cases}
 \end{align}
 which shows that the above equations form a nonhomogeneous system of equations for $W$ and $Z$.
\section{Formulation of the main problem and boundary data}
 We now formulate the problem mathematically in detail by mimicking the real setting of the airfoil problem. Let $\widehat{PE}: r=\varphi(x), x\in [x_1, x_2]$, be a smooth curve and $(\hat{u}(x), \hat{v}(x)), x\in [x_1, x_2]$, is a given velocity distribution on $\widehat{PE}$. Then we define our problem as follows
\subsection{Main Problem}
\textit{Let $\widehat{PE}:r=\varphi(x)~(x\in [x_1, x_2])$ be a smooth streamline of the three-dimensional axisymmetric steady relativistic flow such that it is locally increasing and concave near the point $P$ and $(\hat{u}(x), \hat{v}(x))$ is a given velocity distribution on $\widehat{PE}$ such that $M(x)>1~\forall x\in (x_1, x_3]$ for some $x_3\in (x_1,x_2]$ and $M(x_1)=1$. Then find a smooth sonic curve $\widehat{PD}$ and build a smooth supersonic solution to system \eqref{eq: 2.9} in the angular region of $P$ bounded by $\widehat{PE}$ and $\widehat{PD}$ ; see Figure \ref{fig: 2}.} 
\subsection{Reformulated Problem in terms of angle variables}
We can actually reformulate our main problem in terms of angle variables $(\theta, \varpi)$ as follows.
From Bernoulli's law \eqref{eq: 2.7} and the fact that $\dfrac{da}{dq}<0$, it is easy to see that $a=a(q(x))=a(\hat{u}^2(x)+\hat{v}^2(x))$. Then from $\tan \theta=\dfrac{v}{u}$ and $\sin \omega=\dfrac{1}{M}$, we obtain the data for angle variables $(\theta, \varpi)$ on $\widehat{PE}$ as
 \begin{align}\label{eq: 2.29}
     \theta(x, \varphi(x))&=\tan^{-1}\left(\dfrac{\hat{v}(x)}{\hat{u}(x)}\right)=:\hat{\theta}(x),~~~\varpi(x, \varphi(x))&=\dfrac{\hat{a}(\hat{u}^2(x)+\hat{v}^2(x))\sqrt{1-\hat{a}^2(\hat{u}^2(x)+\hat{v}^2(x))}}{(\hat{u}^2(x)+\hat{v}^2(x))\sqrt{1-(\hat{u}^2(x)+\hat{v}^2(x))}}=:\hat{\varpi}(x).
 \end{align}
 Then we reformulate our problem in terms of angle variables $(\theta, \varpi)$ as:
 \textit{Let us consider a locally increasing smooth streamline $\widehat{PE}:r=\varphi(x)(x\in [x_1, x_2])$ of three-dimensional axisymmetric steady relativistic flow satisfying $\varphi''(x)<0$ in a neighbourhood of $x=x_1$ along which the angle variable $\varpi$ decreases and assign the boundary data $(\theta, \varpi)=(\hat{\theta}, \hat{\varpi})(x)$ on $\widehat{PE}$ such that $\hat{\theta}(x)=\tan^{-1}\varphi'(x), ~~\hat{\varpi}(x)\in(0, 1)~\forall x\in (x_1, x_2]$ and $\hat{\varpi}(x_1)=1$. Then find a smooth sonic curve $\widehat{PD}$ and build a smooth supersonic solution to system \eqref{eq: 2.21} in the angular region of $P$ bounded by $\widehat{PE}$ and $\widehat{PD}$; see Figure \ref{fig: 2}.} 

In order to solve this problem, we assume that the functions $\varphi(x)$ and $\hat{\varpi}(x)$ satisfy \cite{hu2022supersonic}
\begin{align}\label{eq: 2.30}
\begin{cases}
    \varphi(x)\in C^3[x_1, x_2], \hat{\varpi}(x)\in C^2[x_1, x_2],\\
    \varphi(x_1)>0,~~\varphi'(x_1)>0,~~\varphi''(x_1)<0,~~\hat{\varpi}'(x_1)<0,
    \end{cases}
\end{align}
which implies that the curve $r=\varphi(x)$ is increasing and concave  while the angle variable $\varpi$ corresponding to the Mach number decreases near the point $P$. One may note that these assumptions are consistent with the real airfoil setting as well. Since we are focused to develop the existence and regularity of solutions near point o$P$ only, therefore, one may assume without loss of generality
\begin{align}\label{eq: 2.31a}
\begin{cases}
     \varphi(x)\in C^3[x_1, x_2], \hat{\varpi}(x)\in C^2[x_1, x_2],\\
    \varphi_0\leq \varphi(x), ~~\varphi'(x)\leq \varphi_1, \varphi''(x)<0, \hat{\varpi}'(x)<0 ~\forall ~x\in [x_1, x_3]
\end{cases}
\end{align}
for some $x_3\in (x_1, x_2]$, where $\varphi_0$ and $\varphi_1$ are some positive constants. We further assume that $\varphi(x)$ and $\hat{\varpi}(x)$ satisfy
\begin{align}\label{eq: 2.31}
    \left(\dfrac{\varphi''}{1+(\varphi')^2}-\dfrac{4a^2(\varpi)\sqrt{1-(\varpi)^2}}{F_1(\varpi)}(\varpi)'\right)(x)<0 ~~~\forall x\in [x_1, x_3],
\end{align}
which is obviously true near the sonic point $P$. For future use we denote the point $(x_3, \varphi(x_3))$ by $R$ which lie on the curve
$\widehat{PE}$.
\subsection{Boundary data for $W$ and $Z$}
The strategy of this article is to solve system \eqref{eq: 2.28} for $(W, Z)$ in a partial hodograph plane and then return back to the solution via an inverse transformation. Therefore we need to derive the boundary data for $(W, Z)$ on the arc $\widehat{PR}$ using the functions $(\hat{\theta}, \hat{\varpi})(x)$.

Now from \eqref{eq: 2.23} and noting that $\Tilde{\partial}_0=\dfrac{\Tilde{\partial}_++\Tilde{\partial}_-}{2\cos \omega}$, we have
\begin{align}\label{eq: 2.37 a}
    W+Z=-\dfrac{\Tilde{\partial}_0\theta}{\varpi}.
\end{align}
Similarly from \eqref{eq: 2.24}, we have 
\begin{align}\label{eq: 2.38 a}
    W-Z=\dfrac{4a^2(\varpi)\sqrt{1-\varpi^2}\Tilde{\partial}_0\varpi}{\varpi F_1(\varpi)},
\end{align}
which together with \eqref{eq: 2.37 a} implies
\begin{align}\label{eq: 2.32}
    W=\dfrac{2a^2(\varpi)\sqrt{1-\varpi^2}}{\varpi F_1(\varpi)}\Tilde{\partial}_0\varpi-\dfrac{1}{2\varpi}\Tilde{\partial}_0 \theta, ~~~~Z=-\dfrac{2a^2(\varpi)\sqrt{1-\varpi^2}}{\varpi F_1(\varpi)}\Tilde{\partial}_0\varpi-\dfrac{1}{2\varpi}\Tilde{\partial}_0 \theta.
\end{align}
Recalling that the curve $\widehat{PR}$ is a streamline, we have
\begin{align*}
    \Tilde{\partial}_0\theta|_{\widehat{PR}}&=\varphi(x)\cos \hat{\theta}(x) \hat{\theta}'(x)=\varphi(x)\dfrac{\cos \hat{\theta}(x)\varphi''(x)}{1+(\varphi'(x))^2},\\
    \Tilde{\partial}_0\varpi|_{\widehat{PR}}&=\varphi(x)\cos \hat{\theta}(x)\hat{\varpi}'(x),
\end{align*}
which combined with \eqref{eq: 2.32} yields
\begin{align}\label{eq: 2.35}
    \begin{cases}
    W|_{\widehat{PR}}=\dfrac{\varphi(x)\cos \hat{\theta}(x)}{2\hat{\varpi}}\bigg[\dfrac{4a^2(\hat{\varpi})\sqrt{1-\hat{\varpi}^2}}{ F_1(\hat{\varpi})}\hat{\hat{\varpi}}'(x)-\dfrac{\varphi''(x)}{1+(\varphi'(x))^2}\bigg]:=\hat{b}(x),\\
    Z|_{\widehat{PR}}=-\dfrac{\varphi(x)\cos \hat{\theta}(x)}{2\hat{\varpi}}\bigg[\dfrac{4a^2(\hat{\varpi})\sqrt{1-\hat{\varpi}^2}}{ F_1(\hat{\varpi})}\hat{\varpi}'(x)+\dfrac{\varphi''(x)}{1+(\varphi'(x))^2}\bigg]:=\hat{c}(x).
    \end{cases}
\end{align}
For later use, we give the boundary data of $L=\Tilde{\partial}_0 I$ on $\widehat{PR}$
\begin{align}\label{eq: 2.36}
    \Tilde{\partial}_0 I|_{\widehat{PR}}=\dfrac{\varphi(x)\cos \hat{\theta}(x)}{2\hat{\varpi}}\bigg[\dfrac{4a^2(\hat{\varpi})}{F_1(\hat{\varpi})}\hat{\varpi}'(x)\bigg]:=-\hat{d}(x).
\end{align}
Moreover, it suggests by the conditions \eqref{eq: 2.31a} and \eqref{eq: 2.31} that
\begin{align}\label{eq: 2.37}
    \begin{cases}
    \hat{b}(x), ~\hat{c}(x), ~\hat{d}(x)\in C^1([x_1, x_3]),\\
    \hat{m}_0\leq \hat{b}(x),~\hat{c}(x),~ \hat{d}(x)\leq \hat{M}_0, ~\forall ~x\in [x_1, x_3]
    \end{cases}
\end{align}
for some positive constants $\hat{m}_0$ and $\hat{M}_0$.
\section{Existence and regularity of solution in partial hodograph plane}
In this section, we solve the singular system \eqref{eq: 2.28} with the boundary data \eqref{eq: 2.35} under the conditions \eqref{eq: 2.37} near the point $P$ by introducing a partial hodograph transformation.
\subsection{Reformulated problem in partial hodograph plane}
We first reformulate the problem into a new problem by introducing a partial hodograph transformation. We introduce the coordinate transformation $(x, r)\longrightarrow (t, \psi)$ such that
\begin{align}\label{eq: 3.38}
    t=\cos \omega(x, r),~~ \psi=\phi(x, r)-\phi_1,
\end{align}
where $\phi$ is the potential of irrotational relativistic flow such that $\phi_x=u$ and $\phi_r=v$ with $\phi_1=\phi(x_1)$. 

From the transformation we can observe that $\varpi=\sqrt{1-t^2}$. Therefore, we define $F_1(\varpi):=\hat{F}_1(t)$ where
\begin{align}
    \hat{F}_1(t)=4a^2(t)(1-t^2)+2i(t)p''(\rho(t))\gamma_{a(t)}^2f(a(t))>0.
\end{align}
Then by using \eqref{eq: 2.24}, we see that
\begin{align}\label{eq: 3.40}
    J&:=\dfrac{\partial(t, \psi)}{\partial(x, r)}=\begin{vmatrix}
    \dfrac{\partial t}{\partial x}& \dfrac{\partial t}{\partial r}\vspace{0.2 cm}\\
    \dfrac{\partial \psi}{\partial x}& \dfrac{\partial \psi}{\partial r}
    \end{vmatrix}\nonumber\\
    &=\dfrac{a\gamma_a \cos \omega}{\gamma r\sin 2\omega}\dfrac{(\Tilde{\partial}_+\varpi-\Tilde{\partial}_-\varpi)}{\cos \omega}\nonumber\\
    J&=\dfrac{\gamma_a \hat{F}_1(t)(W+Z)}{4a(t)\gamma(t)rt}\neq 0~~\mathrm{for}~0\leq t<1.
\end{align}
We next derive the boundary data of $\phi$ on $\widehat{PR}$. Noting the definitions of $\Tilde{\partial}_i~(i=0, \pm)$ it is easy to obtain that
\begin{align}\label{eq: 4.4}
    \begin{cases}
        \Tilde{\partial}_0 \phi=r \cos \theta. \dfrac{a\gamma_a \cos \theta}{\gamma \sin \omega}+r \sin \theta. \dfrac{a\gamma_a \sin \theta}{\gamma \sin \omega}= \dfrac{a\gamma_a r}{\gamma \sin \omega},\\
        \Tilde{\partial}_+ \phi=r \cos \alpha. \dfrac{a\gamma_a \cos \theta}{\gamma \sin \omega}+r \sin \alpha. \dfrac{a\gamma_a \sin \theta}{\gamma \sin \omega}= \dfrac{a\gamma_a r \cos \omega}{\gamma\sin \omega},\\
        \Tilde{\partial}_- \phi=r \cos \beta. \dfrac{a\gamma_a \cos \theta}{\gamma \sin \omega}+r \sin \beta. \dfrac{a\gamma_a \sin \theta}{\gamma \sin \omega}= \dfrac{a\gamma_a r \cos \omega}{\gamma\sin \omega}.
    \end{cases}
\end{align}
Then noting that the curve $\widehat{PR}$ is a streamline and expression of $\Tilde{\partial}_0 \phi$ in \eqref{eq: 4.4}, one can obtain that
\begin{align}\label{eq: 3.44c}
    \hat{\phi}'(x)=\dfrac{\hat{a}(x)\hat{\gamma_a}(x)}{\hat{\gamma}(x)\hat{\varpi}(x)\cos \hat{\theta}(x)}=\dfrac{\hat{a}(x)\hat{\gamma_a}(x)\sqrt{1+\varphi'(x)^2}}{\hat{\gamma}(x)\hat{\varpi}(x)}>0~\forall~x\in [x_1, x_3],
\end{align}
where $\hat{\phi}(x)=\phi(x, \varphi(x))$. Then we obtain the boundary data of $\hat{\phi}(x)$ on $\widehat{PR}$ as
\begin{align}
    \phi|_{\widehat{PR}}=\hat{\phi}(x)=\phi_1+\displaystyle \int_{x_1}^{x} \dfrac{\hat{a}(s)\hat{\gamma_a}(s)\sqrt{1+\varphi'(s)^2}}{\hat{\gamma}(s)\hat{\varpi}(s)}ds~\forall~x\in [x_1, x_3].
\end{align}

Now using the conditions that $\hat{\phi}'>0$ and $\hat{\varpi}'<0$ it is easy to see that the curve $\widehat{PR}: r=\varphi(x)$ in $x-r$ plane is transformed into a curve $\widehat{P'R'}: \psi=\Tilde{\psi}(t)(t\in [0, t_0])$ in the half plane of $t\geq 0$ defined through a parameter $x$:
\begin{align}\label{eq: 3.44d}
t=\cos \hat{\omega}(x),~~\psi=\hat{\phi}(x)-\phi_1,~~(x\in [x_1, x_3])
\end{align}
such that the number $t_0=\cos \hat{\omega}(x_3)$ is a positive constant. Moreover, since the function $\psi=\hat{\phi}(x)-\phi_1$ is strictly increasing, there exists an inverse function  $x$ such that $x=\hat{x}(\psi)(\psi\in [0, \psi_0])$, where $\psi_0=\hat{\phi}(x_3)-\phi_1$. Now we denote $\hat{b}(\psi)=\hat{b}(x(\psi)), ~\hat{c}(\psi)=\hat{c}(x(\psi)),~\hat{d}(\psi)=\hat{d}(x(\psi))$. Then we obtain the boundary data of $(W, Z, L)$ on $\widehat{P'R'}$ such that
\begin{align}\label{eq: 2.42}
    W|_{\widehat{P'R'}}=\hat{b}(\psi),~~Z|_{\widehat{P'R'}}=\hat{c}(\psi),~~L|_{\widehat{P'R'}}=-\hat{d}(\psi)~~ \forall~ \psi\in [0, \psi_0].
\end{align}
Then it is straightforward to see from \eqref{eq: 2.37} that
\begin{align}\label{eq: 2.43}
\begin{cases}
    \hat{b}(\psi), \hat{c}(\psi), \hat{d}(\psi)\in C^1[0, \psi_0],\\
    \hat{m}_0\leq \hat{b}(\psi),~\hat{c}(\psi), ~\hat{d}(\psi) \leq \hat{M}_0,~~ \forall ~\psi\in [0, \psi_0]
  \end{cases}  
\end{align}
for some positive constants $\hat{m}_0$ and $\hat{M}_0$.

Further, in these new coordinates, the weighted normalized derivatives now become
\begin{align}\label{eq: 3.39}
    \Tilde{\partial}_+&=-\dfrac{2F(t)W}{t}\bigg[\dfrac{\partial}{\partial t}-\dfrac{a r \gamma_a t^2}{2\gamma F(t)W  \sqrt{1-t^2}}\dfrac{\partial}{\partial \psi}\bigg],~~~\Tilde{\partial}_-&=\dfrac{2F(t)Z}{t}\bigg[\dfrac{\partial}{\partial t}+\dfrac{a r \gamma_a t^2}{2\gamma F(t)Z  \sqrt{1-t^2}}\dfrac{\partial}{\partial \psi}\bigg]
\end{align}
where $F(t)=\dfrac{(1-t^2)\hat{F}_1(t)}{4\hat{a}^2(t)}>0.$

Then we derive the characteristic decompositions of $(W, Z)$ in the partial hodograph $(t, \psi)$ plane. Using the normalized derivatives \eqref{eq: 3.39} in $t-\psi$ plane in \eqref{eq: 2.28}, it is easy to obtain the reformulated characteristic decompositions of $(W, Z)$ of the form
\begin{align}\label{eq: 3.44}
\begin{cases}
    &\hspace{-0.4 cm}W_t+\dfrac{a r \gamma_a t^2}{2\gamma F(t)Z  \sqrt{1-t^2}}W_\psi=\bigg[1+\dfrac{i\gamma_a^2p''(\rho)f(a)}{2a^2}\bigg]\dfrac{W}{ZF(t)}\left(\dfrac{W-Z}{2t}\right)\\
    &\hspace{2 cm}+\bigg[\dfrac{i\gamma_a^2 p''(\rho)f(a)}{2a^2}+2(1-t^2)+\dfrac{\sin \beta+4 \sin \theta t}{4Z}\bigg]\dfrac{Wt}{F(t)}
    -\dfrac{\sin \beta}{4F(t)Z}Zt+\dfrac{\sin^2 \theta}{2F(t)Z}t,\\
    &\hspace{-0.4 cm}Z_t-\dfrac{a r \gamma_a t^2}{2\gamma F(t)W  \sqrt{1-t^2}}Z_\psi=\bigg[1+\dfrac{i\gamma_a^2p''(\rho)f(a)}{2a^2}\bigg]\dfrac{Z}{WF(t)}\left(\dfrac{Z-W}{2t}\right)\\
    &\hspace{2 cm}+\bigg[\dfrac{i\gamma_a^2 p''(\rho)f(a)}{2a^2}+2(1-t^2)-\dfrac{\sin \alpha+4 \sin \theta t}{4W}\bigg]\dfrac{Zt}{F(t)}+\dfrac{\sin \alpha}{4F(t)W}Wt+\dfrac{\sin^2 \theta}{2F(t)W}t,\\
\end{cases}    
\end{align}
It is easy to see that \eqref{eq: 3.44} is not a closed system because there are two unknown functions $r(t, \psi)$ and $\theta(t, \psi)$ in the system. In order to close the system, we need to add the characteristic equations and boundary data of $r$ and $\theta$ to the system. First of all, from the definitions of weighted directional derivatives, we have the characteristic decompositions of the form:
\begin{align*}
    \Tilde{\partial}_+r=r\sin \alpha,~~\Tilde{\partial}_-r=r\sin \beta,
\end{align*}
from which one can have
\begin{align}\label{eq: 3.44a}
    \partial^- r=-\dfrac{tr \sin \alpha}{2FW},~~\partial^+ r=\dfrac{tr \sin \beta}{2FZ},
\end{align}
where $\partial^{\pm}=\partial_t+\lambda_{\pm}\partial_\psi$ such that $\lambda_+= \dfrac{a r \gamma_a t^2}{2\gamma F(t)Z  \sqrt{1-t^2}}$ and $\lambda_-= -\dfrac{a r \gamma_a t^2}{2\gamma F(t)W  \sqrt{1-t^2}}$.

Similarly, from \eqref{eq: 2.23} we get the characteristic decomposition of $\theta$ of the form
\begin{align}\label{eq: 3.44b}
    \partial^- \theta=\dfrac{t\sqrt{1-t^2} (\sin \theta+2tW)}{2FW},~~\partial^+ \theta=\dfrac{t\sqrt{1-t^2} (\sin \theta-2tZ)}{2FZ}.
\end{align}
Also, the boundary data of $r$ and $\theta$ on $\widehat{P'R'}$ are
\begin{align}\label{eq: 3.56}
    r|_{\widehat{P'R'}}=\varphi(\hat{x}(\psi))=: \hat{\varphi}(\psi), ~~\theta|_{\widehat{P'R'}}=\hat{\theta}(\hat{x}(\psi))=:\hat{\theta}(\psi),
\end{align}
which satisfy
\begin{align}\label{eq: 3.57}
    \begin{cases}
        \hat{\varphi}(\psi), \hat{\theta}(\psi)\in C^1([0, \psi_0]),\\
        \varphi_0\leq \hat{\varphi}(\psi)\leq \varphi_1 \forall~\psi\in [0, \psi_0].
    \end{cases}
\end{align}
Combining \eqref{eq: 3.44}, \eqref{eq: 3.44a} and \eqref{eq: 3.44b}, we obtain a closed system of $(W,Z, r, \theta)(t, \psi)$ of the form
\begin{align}\label{eq: 3.45}
\begin{cases}
    &\hspace{-0.4 cm}\partial^+W=\bigg[1+\dfrac{i\gamma_a^2p''(\rho)f(a)}{2a^2}\bigg]\dfrac{W}{ZF(t)}\left(\dfrac{W-Z}{2t}\right)+\bigg[\dfrac{i\gamma_a^2 p''(\rho)f(a)}{2a^2}+2(1-t^2)+\dfrac{\sin \beta+4 \sin \theta t}{4Z}\bigg]\dfrac{Wt}{F(t)}\\
    &\hspace{4 cm}-\dfrac{\sin \beta}{4F(t)Z}Zt+\dfrac{t\sin^2 \theta}{2F(t)Z},\\
    &\hspace{-0.4 cm}\partial^-Z=\bigg[1+\dfrac{i\gamma_a^2p''(\rho)f(a)}{2a^2}\bigg]\dfrac{Z}{WF(t)}\left(\dfrac{Z-W}{2t}\right)+\bigg[\dfrac{i\gamma_a^2 p''(\rho)f(a)}{2a^2}+2(1-t^2)-\dfrac{\sin \alpha+4 \sin \theta t}{4W}\bigg]\dfrac{Zt}{F(t)}\\
    &\hspace{4 cm}+\dfrac{\sin \alpha}{4F(t)W}Wt+\dfrac{t\sin^2 \theta}{2F(t)W},\\
    &\hspace{-0.4 cm}\partial^- r=-\dfrac{tr \sin \alpha}{2FW},\vspace{0.2 cm}\\
    &\hspace{-0.4 cm}\partial^- \theta=\dfrac{t\sqrt{1-t^2} (\sin \theta+2tW)}{2FW}
 \end{cases}   
\end{align}
with the boundary data \eqref{eq: 2.42} and \eqref{eq: 3.56}. It is easy to observe that the two eigenvalues of the system \eqref{eq: 3.45} are $\lambda_\pm$ defined as before. 

In the upcoming subsections, we try to obtain the existence and regularity of a solution $(W, Z, r, \theta)$ of the closed system $\eqref{eq: 3.45}$ with the boundary data \eqref{eq: 2.42} and \eqref{eq: 3.56} in the neighbourhood $P'F'D'$ of the point $P'$, where $F'$ is a point $(\bar{\delta}, \Tilde{\psi}(\bar{\delta}))$ lying on the streamline $\widehat{P'R'}$ such that $\bar{\delta}\in (0, t_0]$ is a small positive number; see Figure \ref{fig: region}.
\subsection{A strong determinate domain and a priori estimates}
 In this subsection, we construct a strong determinate domain $\Omega$ for system \eqref{eq: 3.45}, which is not easy due to the nonlinearity of the system.
 
Noting the fact that $0<a^2<1,~ p''(\rho)>0,~i(\rho)>0,~0<\gamma_a<1,~f(a)>0$, let us set
\begin{align}\label{eq: 3.59}
    \overline{K}=1+&\max\bigg\{1+\left(\dfrac{2\hat{m}_0+1}{\hat{m}_0}\right)\bigg[\dfrac{2a^2}{i\gamma_a^2 p''(\rho)f(a)}\bigg], 1+\left(\dfrac{2\hat{m}_0^2+2\hat{m}_0+2}{\hat{m}_0^2}\right)\bigg[\dfrac{2a^2}{i\gamma_a^2p''(\rho)f(a)}\bigg]\bigg\}< \infty.
\end{align} 
Further, we choose $\delta_0=\min\{1/\sqrt{2}, t_0\}$ such that
\begin{align}
    F(t)=(1-t^2)\bigg[\dfrac{i\gamma_a^2p''(\rho)f(a)}{2a^2}+1-t^2\bigg]\geq \dfrac{i\gamma_a^2p''(\rho)f(a)}{2a^2}\geq\dfrac{\kappa}{2}>0 ~\forall~t\in (0, \delta_0]
\end{align}
for some $\kappa>0$. Moreover, we choose a positive number $\delta_1 \leq \delta_0$ such that
\begin{align}
    e^{\overline{K}\delta_1^2}\leq 2,~~\delta_1\leq \dfrac{\kappa \hat{m}_0}{3}.
\end{align}
We next derive the slope of curve $\widehat{P'R'}\cap \{t\leq \delta_0\}$. From \eqref{eq: 3.44c} and \eqref{eq: 3.44d} we have
\begin{align}\label{eq: 3.61}
    \Tilde{\psi}'(t)=-\hat{\phi}'(x).\dfrac{\cos \hat{\omega}}{\hat{\varpi}(x)\hat{\varpi}'(x)}=-\dfrac{\hat{a}\hat{\gamma}_a\sqrt{1+\varphi'^2}}{\hat{\gamma} \hat{\varpi}'\hat{\varpi}^2}t.
\end{align}
Then if we denote $\Tilde{m}=\underset{-\psi\in [0, \tilde{\psi}(\delta_0)]}\min \dfrac{\hat{a}\hat{\gamma}_a\sqrt{1+\varphi'^2}}{\hat{\gamma} \hat{\varpi}'\hat{\varpi}^2}(\hat{x}(\psi))>0$ and set\\ $\widetilde{K}=\max\bigg\{\dfrac{8\varphi_1}{i\gamma_a^2 p''(\rho)f(a) \hat{m}_0}exp\left(\dfrac{2t_0^2}{i\gamma_a^2 p''(\rho)f(a) \hat{m}_0}\right)\bigg\}< \infty$ such that there exists a $\delta_2\leq \delta_0$ which satisfy $\widetilde{K}\delta_2<\Tilde{m}$.

Let $\delta=\min\{\delta_1, \delta_2\}$. Then we consider the curve $\psi=\bar{\psi}(t)$ defined by 
\begin{align}
    \bar{\psi}(t)=\tilde{\psi}(\delta)-\dfrac{\widetilde{K}}{3}\delta^3+\dfrac{\widetilde{K}}{3}t^3 ~\forall~t\in [0, \delta]. 
\end{align}
Then by integrating $\eqref{eq: 3.61}$, it is easy to see that
\begin{align}
    \tilde{\psi}(\delta)\geq \displaystyle \int_{0}^{\delta} \widetilde{K}t dt=\dfrac{\Tilde{m}}{2}\delta^2>\dfrac{\widetilde{K}\delta_2}{2}\delta^2>\dfrac{\widetilde{K}}{3}\delta^3.
\end{align}
so that $\bar{\psi}(0)>0$.
\begin{figure}
    \centering
    \includegraphics[width= 3.5 in]{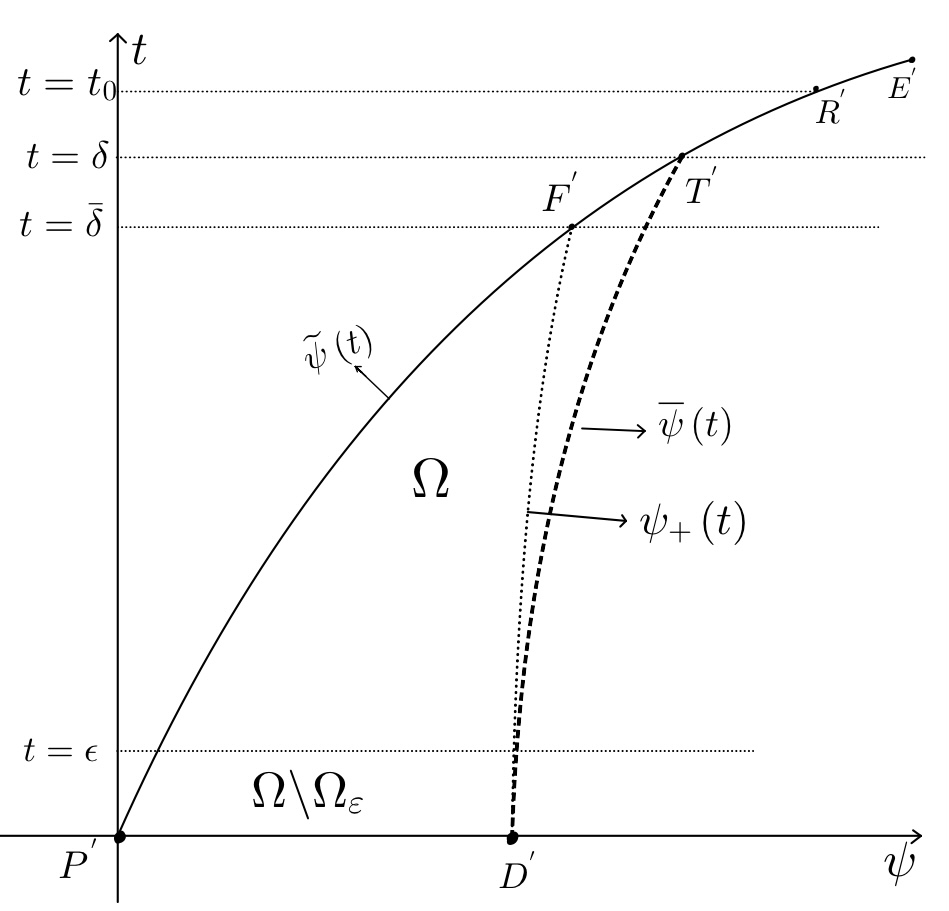}
    \caption{The strong determinate domain $\Omega$}
    \label{fig: region}
\end{figure}
Let us denote the point $(0, \bar{\psi}(0))$ by $D'$ and the point $(\delta, \tilde{\psi}(\delta))$ by $T'$. Further, let $\Omega$ be the domain bounded by the curves $\widehat{P'T'}$, $\widehat{D'T'}$ and the degenerate line $\widehat{P'D'}$. Moreover, $\epsilon\in (0, \delta]$ be an arbitrary constant such that we denote $\Omega_\epsilon=\Omega\cap \{(t, \psi)|t\geq \epsilon\}$; see Figure \ref{fig: region}. Then we have the following Lemma:
\subsubsection{{Upper and lower bounds of $W, Z$ and $r$ in $t-\psi$ plane}}
\begin{l1}\label{l: 4.1}
Let us assume that the conditions \eqref{eq: 2.43} and \eqref{eq: 3.57} are satisfied such that there exists a $C^1$ solution $(W, Z, r, \theta)(t, \psi)$ of the system \eqref{eq: 3.45} with the boundary data \eqref{eq: 2.42} and \eqref{eq: 3.56} in the domain $\Omega_\epsilon$. Then
\begin{align}
    \begin{cases}
        \dfrac{\hat{m}_0}{2}< W(t, \psi), Z(t, \psi)< \hat{M}_0+1,\\
        \varphi_0exp \left(-\dfrac{t_0^2}{\kappa \hat{m}_0}\right)\leq r(t, \psi)\leq \varphi_1 exp \left(\dfrac{t_0^2}{\kappa \hat{m}_0}\right) 
    \end{cases}
\end{align}
for all $(t, \psi)\in \Omega_\epsilon$.
\end{l1}
\begin{proof}
In order to prove this Lemma, we need to prove that the region $\mathbb{I}:=\left(\dfrac{\hat{m}_0}{2}, \hat{M}_0+1\right)\times \left(\dfrac{\hat{m}_0}{2}, \hat{M}_0+1\right)$ is an invariant region for $(W, Z)(t, \psi)$.

Let us consider the level curve $t=\epsilon'(\epsilon'\in [\epsilon, \delta])$ and move it from $t=\delta$ to $t=\epsilon$. We assume, on the contrary, that the region $\mathbb{I}$ is not an invariant region and consider that the point $A$ is the first time on one such level set $t=\epsilon'$ such that one of $W$ and $Z$ touches the lower boundary of $\mathbb{I}$. Without loss of generality, we assume that $W|_A=\dfrac{\hat{m}_0}{2}$ and $(W, Z)\in \mathbb{I}$ for any $(t, \psi)\in \Omega_\epsilon\cap \{t>\epsilon'\}$. Now it is easy to check by applying the third equation of \eqref{eq: 3.45} that 
\begin{align*}
    \varphi_0exp \left(-\dfrac{t_0^2}{\kappa \hat{m}_0}\right)\leq r(t, \psi)\leq \varphi_1 exp \left(\dfrac{t_0^2}{\kappa \hat{m}_0}\right) ~~\forall ~(t, \psi)\in \Omega_\epsilon\cap \{t>\epsilon'\}.
\end{align*}
Hence, we draw a positive characteristic $l_1^+$ from the point $A$ up to a point $A_1$ on the boundary curve $\widehat{P'T'}$. Further, we set
\begin{align}
    \widetilde{W}=W-g(t), ~~\widetilde{Z}=Z-g(t),
\end{align}
where $~~g(t)=\dfrac{\hat{m}_0}{2}e^{\overline{K}t^2}.$

By the choice of $\delta_1$, we have $\widetilde{W}|_{A_1}>0.$ Also, it follows by $W|_{A}=\hat{m}_0/2$ that $\widetilde{W}|_{A}<0$, which implies by the continuity of $W$ and $g$ that there exists a point $B$ lying between $A$ and $A_1$ on $l_1^+$ such that $\widetilde{W}|_{B}=0$ and $\widetilde{W}|_{\widehat{BA_1}}\geq 0$. Therefore, $\partial^+ \widetilde{W}|_{B}\geq 0$.

However, using the equations of $W$ and $Z$ from the system \eqref{eq: 3.45} and performing a direct calculation, one can obtain the following characteristic decompositions of $(\widetilde{W}, \widetilde{Z})$ 
\begin{align}\label{eq: 3.67}
    \begin{cases}
&\hspace{-0.5 cm}\partial^+\widetilde{W}=\bigg[1+\dfrac{i\gamma_a^2p''(\rho)f(a)}{2a^2}\bigg]\dfrac{(\widetilde{W}+g)}{(\widetilde{Z}+g)F(t)}\left(\dfrac{\widetilde{W}-\widetilde{Z}}{2t}\right)+\bigg[\dfrac{i\gamma_a^2 p''(\rho)f(a)}{2a^2}+2(1-t^2)+\dfrac{\sin \beta+4 \sin \theta t}{4(\widetilde{Z}+g)}\bigg]\dfrac{\widetilde{W}t}{F(t)}\\
    &\hspace{3 cm}+\dfrac{\widetilde{Z}t}{F(t)(\widetilde{Z}+g)}\Upsilon_1+\dfrac{t}{F(t)(\widetilde{Z}+g)} \Upsilon_2,\vspace{0.3 cm}\\
        &\hspace{-0.5 cm}\partial^-\widetilde{Z}=\bigg[1+\dfrac{i\gamma_a^2p''(\rho)f(a)}{2a^2}\bigg]\dfrac{(\widetilde{Z}+g)}{(\widetilde{W}+g)F(t)}\left(\dfrac{\widetilde{Z}-\widetilde{W}}{2t}\right)+\bigg[\dfrac{i\gamma_a^2 p''(\rho)f(a)}{2a^2}+2(1-t^2)-\dfrac{\sin \alpha+4 \sin \theta t}{4(\widetilde{W}+g)}\bigg]\dfrac{\widetilde{Z}t}{F(t)}\\
        &\hspace{3 cm}+\dfrac{\widetilde{W}t}{F(t)(\widetilde{W}+g)}\Upsilon_3+\dfrac{t}{F(t)(\widetilde{W}+g)} \Upsilon_4,\\
    \end{cases}
\end{align}
where 
\begin{align*}
    \Upsilon_1&=\bigg[\dfrac{i\gamma_a^2 p''(\rho)f(a)}{2a^2}+2(1-t^2)\bigg]g(t)-2F(t)g(t)\overline{K}-\dfrac{\sin \beta}{4},\\
    \Upsilon_2&=\bigg[\dfrac{i\gamma_a^2 p''(\rho)f(a)}{2a^2}+2(1-t^2)\bigg]g^2(t)-2F(t)g^2(t)\overline{K}+g(t)\sin \theta t+\dfrac{\sin^2 \theta}{2},\\
    \Upsilon_3&=\bigg[\dfrac{i\gamma_a^2 p''(\rho)f(a)}{2a^2}+2(1-t^2)\bigg]g(t)-2F(t)g(t)\overline{K}+\dfrac{\sin \alpha}{4},\\
    \Upsilon_4&=\bigg[\dfrac{i\gamma_a^2 p''(\rho)f(a)}{2a^2}+2(1-t^2)\bigg]g^2(t)-2F(t)g^2(t)\overline{K}-g(t)\sin \theta t+\dfrac{\sin^2 \theta}{2}.
\end{align*}
For $\Upsilon_1$ and $\Upsilon_3$, one has
\begin{align*}
    \Upsilon_1, ~\Upsilon_3&\leq \bigg[\dfrac{i\gamma_a^2 p''(\rho)f(a)}{2a^2}+2(1-t^2)\bigg]g(t)-2F(t)g(t)\overline{K}+\dfrac{1}{4},\\
    & \leq g\bigg[\dfrac{i\gamma_a^2 p''(\rho)f(a)}{2a^2}+2+\dfrac{1}{4g}-2F(t)\overline{K}\bigg]\\
    &\leq g(t)\bigg[\dfrac{i\gamma_a^2 p''(\rho)f(a)}{2a^2}+2+\dfrac{1}{\hat{m}_0}-\dfrac{i\gamma_a^2 p''(\rho)f(a)}{2a^2}\overline{K}\bigg]\\
    &=\dfrac{i\gamma_a^2 p''(\rho)f(a)}{2a^2}g(t)\left(1+\left(\dfrac{2\hat{m}_0+1}{\hat{m}_0}\right)\bigg[\dfrac{2a^2}{i\gamma_a^2p''(\rho)f(a)}\bigg]-\overline{K}\right)<0
\end{align*}
by the choice of $\overline{K}$ in \eqref{eq: 3.59}. 

Similarly, for $\Upsilon_2$ and $\Upsilon_4$ also, we have
\begin{align*}
   \Upsilon_2, &\Upsilon_4  \leq \left(\dfrac{i\gamma_a^2 p''(\rho)f(a)}{2a^2}+2\right)g^2-\dfrac{i\gamma_a^2 p''(\rho)f(a)}{2a^2}g^2\overline{K}+g(t)+\dfrac{1}{2}\\
    & \leq \dfrac{i\gamma_a^2 p''(\rho)f(a)}{2a^2}g^2(t)\left(1+\left(\dfrac{2\hat{m}_0^2+2\hat{m}_0+2}{\hat{m}_0^2}\right)\bigg[\dfrac{2a^2}{i\gamma_a^2p''(\rho)f(a)}\bigg]-\overline{K}\right)<0
\end{align*}
by the choice of $\overline{K}$ in \eqref{eq: 3.59}. 

Therefore, using the equation of $\widetilde{W}$ from \eqref{eq: 3.67} and noting the facts that $\Upsilon_1, \Upsilon_2<0$, one can conclude that $\partial^+\widetilde{W}|_{B}<0$, which leads to a contradiction.

Similarly, if there exists a point $B'$ lying on a level curve $t=\epsilon'$ in $\Omega_\epsilon$ such that one of $W$ or $Z$ touches the upper boundary of $\mathbb{I}$. Again, without loss of generality, we assume that $W|_{B'}=\hat{M}_0+1$ and $(W, Z)\in \mathbb{I}$ for any $(t, \psi)\in \Omega_\epsilon\cap\{t>\epsilon'\}$. Thus, we can draw a positive characteristic $l_2^+$ from $B'$ up to a point $B_1$ lying on the boundary curve $\widehat{P'T'}$. On the curve $\widehat{B'B_1}$, we have $W\leq \hat{M}_0+1$, which implies that $\partial^+W|_{B'}\leq 0$. However, using the equation for $W$ from \eqref{eq: 3.45}, one can obtain
\begin{align*}
    \partial^+W|_{B'}\geq \bigg[\dfrac{i\gamma_a^2 p''(\rho)f(a)}{2a^2}+2(1-t^2)+\dfrac{\sin \theta t}{Z}\bigg]\dfrac{Wt}{F(t)}+\dfrac{\sin^2 \theta}{2F(t)Z}t>0,
\end{align*}
where we have used the fact $\bigg[\dfrac{i\gamma_a^2 p''(\rho)f(a)}{2a^2}+2(1-t^2)\bigg]Z+t \sin \theta \geq \dfrac{\kappa}{2}\hat{m}_0-\delta |\sin \theta|>0$, which is true by the choice of $\delta_1$.

The above two conclusions lead to a contradiction that proves that the region $\mathbb{I}$ is an invariant region of $(W, Z)$. Furthermore, the estimates of $r(t, \psi)$ can be easily obtained using the third equation of \eqref{eq: 3.45}. Therefore, the proof of the Lemma is completed.
\end{proof}
Now we consider a function class $\mathcal{S}(\overline{\Omega})$ which incorporates all vector functions $\mathbf{F}=(f_1, f_2, f_3, f_4)^T:\overline{\Omega}\longrightarrow \mathbb{R}^4$ satisfying the following properties:
\begin{align}
    \begin{cases}
        \hspace{-0.1 cm}(P_1):~f_1, f_2\in C^1(\overline{\Omega}/\{t=0\}), f_3, f_4\in C^1(\overline{\Omega});\\
        \hspace{-0.1 cm}(P_2):~(f_1, f_2, f_3, f_4)^T(t, \Tilde{\psi}(t))=(\hat{b}(\Tilde{\psi}(t)), \hat{c}(\Tilde{\psi}(t)), \hat{\varphi}(\Tilde{\psi}(t)), \hat{\theta}(\Tilde{\psi}(t)))^T~\forall t\in[0, \delta];\\
        \hspace{-0.1 cm}(P_3):~\dfrac{\hat{m}_0}{2}< f_1, f_2<\hat{M}_0+1, ~~\varphi_0 exp\left(-\dfrac{t_0^2}{\kappa \hat{m}_0}\right)\leq f_3\leq \varphi_1 exp\left(\dfrac{t_0^2}{\kappa \hat{m}_0}\right).
    \end{cases}
\end{align}
It can be easily observed from Lemma \ref{l: 4.1} that $\mathcal{S}(\overline{\Omega})$ is not empty. Further, based on the expression of $\lambda_+$, we define the curve $\psi=\psi_+(t; \xi, \eta)$ by
\begin{align}
\begin{cases}
    \dfrac{d\psi_+(t; \xi, \eta)}{dt}=\dfrac{ra\gamma_a t^2}{2\gamma FZ\sqrt{1-t^2}}(t; \xi, \eta),\\
    \psi_+(\xi; \xi, \eta)=\eta,
\end{cases}    
\end{align}
for $t\geq \xi$, where $(\xi, \eta)$ is an arbitrary point in $\Omega$. Then we proceed to prove that $\Omega$ is a strong determinate domain in the following Lemma:
\begin{l1}
If the solution $(W, Z, r, \theta)(t, \psi)$ of the system \eqref{eq: 3.45} with the boundary data \eqref{eq: 2.42} and \eqref{eq: 3.56} belongs to the function class $\mathcal{S}(\overline{\Omega})$ then $\Omega$ is a strong determinate domain.
\end{l1}
\begin{proof}
In order to prove this Lemma, it is enough to prove that the curve $\psi=\psi_+(t; \xi, \eta)$ intersects only with the curve $\widehat{P'T'}$. We prove this by proving that the slope of curve $\psi=\psi_+(t)$ is strictly smaller than $\psi=\bar{\psi}(t)$ at any point on $\widehat{D'T'}$. Indeed, by using \eqref{eq: 3.59}
\begin{align}
    \dfrac{ra\gamma_a t^2}{2\gamma FZ\sqrt{1-t^2}}\leq \dfrac{t^2}{2.\dfrac{\kappa}{2}\dfrac{\hat{m}_0}{2}.\dfrac{1}{\sqrt{2}}}\varphi_1 exp\left(\dfrac{t_0^2}{\kappa \hat{m}_0}\right)<\widetilde{K}t^2=\bar{\psi}'(t)
\end{align}
for $t\in (0, \delta)$ which implies that the domain $\Omega$ is a strong determinate domain.
\end{proof}

\subsection{Existence of Solutions}\label{s: 4.3}
To establish the existence of solutions, we first need to derive a priori $C^1$ estimates. For this purpose, we first introduce $\overline{W}=\dfrac{1}{W},~ \overline{Z}=\dfrac{1}{Z}$ to convert the system \eqref{eq: 3.45} into the following form
\begin{align}\label{eq: 3.71}
\begin{cases}
    \bar{\partial}^+\overline{W}=\dfrac{\overline{W}-\overline{Z}}{2t}+tH_1(\overline{W}, \overline{Z}, \theta, t),\vspace{0.2 cm}\\
    \bar{\partial}^-\overline{Z}=\dfrac{\overline{Z}-\overline{W}}{2t}+tH_2(\overline{W}, \overline{Z}, \theta, t),\vspace{0.2 cm}\\
    \bar{\partial}^-r=-\dfrac{\overline{W}tr \sin \alpha}{2F},\vspace{0.2 cm}\\
    \bar{\partial}^-\theta=\dfrac{t\sqrt{1-t^2}(\overline{W} \sin \theta+2t)}{2F},
\end{cases}
\end{align}
where $\bar{\partial}^\pm=\partial_t+\bar{\lambda}_\pm \partial_\psi$ with $\bar{\lambda}_+=\dfrac{\overline{Z}ra\gamma_a}{2\gamma F\sqrt{1-t^2}}t^2$, $\bar{\lambda}_-=-\dfrac{\overline{W}ra\gamma_a}{2\gamma F\sqrt{1-t^2}}t^2$,\\
and 
\begin{align*}
    H_1=&\bigg[\dfrac{i\gamma_a^2 p''(\rho)f(a)}{2a^2}+2-t^2\bigg]\dfrac{(\overline{W}-\overline{Z})}{2F(t)}-\bigg[\dfrac{i\gamma_a^2p''(\rho)f(a)}{2a^2}+2(1-t^2)\bigg]\dfrac{\overline{W}}{F(t)}\\
    &-\dfrac{\overline{W}}{2F}\bigg\{\dfrac{\sin \beta}{2}(\overline{Z}-\overline{W})+2t\overline{Z}\sin \theta+\overline{W}\cdot \overline{Z}\sin^2\theta\bigg\},\\
    H_2=&\bigg[\dfrac{i\gamma_a^2 p''(\rho)f(a)}{2a^2}+2-t^2\bigg]\dfrac{(\overline{Z}-\overline{W})}{2F(t)}-\bigg[\dfrac{i\gamma_a^2p''(\rho)f(a)}{2a^2}+2(1-t^2)\bigg]\dfrac{\overline{Z}}{F(t)}\\
    &-\dfrac{\overline{Z}}{2F}\bigg\{\dfrac{\sin \alpha}{2}(\overline{Z}-\overline{W})-2t\overline{W}\sin \theta+\overline{W} \cdot \overline{Z}\sin^2\theta\bigg\},
\end{align*}
We now use the following commutator relation
\begin{align}\label{eq: 3.72a}
    \bar{\partial}^-\bar{\partial}^+-\bar{\partial}^+\bar{\partial}^-=\dfrac{\bar{\partial}^-\bar{\lambda}_+-\bar{\partial}^+\bar{\lambda}_-}{\bar{\lambda}_+-\bar{\lambda}_-}(\bar{\partial}^+-\bar{\partial}^-)
\end{align}
to obtain the equations of $\bar{\partial}^-\overline{W}$ and $\bar{\partial}^+\overline{Z}$ as follows:
\begin{align}\label{eq: 3.72}
\begin{cases}
    \bar{\partial}^+\bar{\partial}^- \overline{W}=\bar{\partial}^-\bar{\partial}^+ \overline{W}+\dfrac{\bar{\partial}^-\bar{\lambda}_+-\bar{\partial}^+\bar{\lambda}_-}{\bar{\lambda}_+-\bar{\lambda}_-}(\bar{\partial}^-\overline{W}-\bar{\partial}^+\overline{W}),\\
    \bar{\partial}^-\bar{\partial}^+ \overline{Z}=\bar{\partial}^+\bar{\partial}^- \overline{Z}+\dfrac{\bar{\partial}^-\bar{\lambda}_+-\bar{\partial}^+\bar{\lambda}_-}{\bar{\lambda}_+-\bar{\lambda}_-}(\bar{\partial}^+\overline{Z}-\bar{\partial}^-\overline{Z}).
\end{cases}    
\end{align}
A routine calculation now yields
\begin{align}\label{eq: 3.73}
    \dfrac{\bar{\partial}^-\bar{\lambda}_+-\bar{\partial}^+\bar{\lambda}_-}{\bar{\lambda}_+-\bar{\lambda}_-}=\dfrac{2}{t}+th,
\end{align}
where
\begin{align*}
    h=&\bigg[\dfrac{i\gamma_a^2p''(\rho)f(a)}{2a^2}+4(1-t^2)\bigg]\dfrac{1}{F(t)}+\dfrac{1}{1-t^2}+\dfrac{H_1+H_2}{\overline{W}+\overline{Z}}-\dfrac{\overline{W}\cdot\overline{Z}\sqrt{1-t^2}\cos \theta}{\overline{W}+\overline{Z}}\\
    &~~~~~~~~~~~~-\dfrac{i\gamma_a^2t(1-t^2)}{F^2}\left(\gamma_a^2 p''(\rho)+\tau'(\rho)\right),
\end{align*}
where $\tau= \dfrac{i\gamma_a^4 p''(\rho)}{2a^2}$ such that $\tau'$ is bounded in the domain $\Omega$.

Moreover, making use of \eqref{eq: 3.71}, we find that
\begin{align}\label{eq: 3.74}
    \bar{\partial}^-\bar{\partial}^+ \overline{W}=\dfrac{\bar{\partial}^-\overline{W}-\bar{\partial}^-\overline{Z}}{2t}-\dfrac{\overline{W}-\overline{Z}}{2t^2}+H_1+tf_1\bar{\partial}^-\overline{W}+tf_2\bar{\partial}^-\overline{Z}+tf_3\bar{\partial}^-\theta+tf_4,
\end{align}
where
\begin{align*}
    f_1&=\dfrac{\partial H_1}{\partial \overline{W}}=\dfrac{1}{2F}\Bigg\{3t^2-2-\dfrac{i\gamma_a^2p''(\rho)f(a)}{2a^2}-2t\overline{Z}\sin \theta-\dfrac{\sin \beta}{2}\overline{Z}-2 \overline{W}. \overline{Z}\sin^2\theta+\sin \beta \overline{W}\Bigg\},\\
    f_2&=\dfrac{\partial H_1}{\partial \overline{Z}}=\dfrac{1}{2F}\Bigg\{t^2-2-\dfrac{i\gamma_a^2p''(\rho)f(a)}{2a^2}-2t\overline{W}\sin \theta-\dfrac{\sin \beta}{2}\overline{W}- \overline{W}^2\sin^2\theta\Bigg\},\\
    f_3&=\dfrac{\partial H_1}{\partial \theta}=-\dfrac{\overline{W}}{2F}\Bigg\{2t\overline{Z}\cos \theta+\dfrac{\cos \beta}{2}(\overline{Z}- \overline{W})+\overline{W}.\overline{Z}\sin 2\theta \Bigg\},\\
    f_4&=\dfrac{\partial H_1}{\partial t}=\dfrac{H_1}{F}\Bigg\{\bigg[\left(\dfrac{i\gamma_a^2p''(\rho)f(a)}{2a^2}+2(1-t^2)\right)2t-i\gamma_a^2t^2(1-t^2)\left(\gamma_a^2p''(\rho)+\tau'(\rho)\right)\bigg]\\
    &~~~~~~~~~~~~~+3t\overline{W}+t\overline{Z}-\overline{W}.\overline{Z}\sin \theta+i\gamma_a^2t^2\left(\gamma_a^2p''(\rho)+\tau'(\rho)\right)\left(\dfrac{\overline{W}+\overline{Z}}{2}\right)\Bigg\}
\end{align*}
and 
\begin{align}\label{eq: 3.75}
    \bar{\partial}^+\bar{\partial}^- \overline{Z}=\dfrac{\bar{\partial}^+\overline{Z}-\bar{\partial}^+\overline{W}}{2t}-\dfrac{\overline{Z}-\overline{X}}{2t^2}+H_2+tg_1\bar{\partial}^+\overline{Z}+tg_2\bar{\partial}^+\overline{W}+tg_3\bar{\partial}^+\theta+tg_4,
\end{align}
where
\begin{align*}
    g_1&=\dfrac{\partial H_2}{\partial \overline{Z}}=\dfrac{1}{2F}\Bigg\{3t^2-2-\dfrac{i\gamma_a^2p''(\rho)f(a)}{2a^2}+2t\overline{W}\sin \theta+\dfrac{\sin \alpha}{2}\overline{W}-2 \overline{W}. \overline{Z}\sin^2\theta-\sin \alpha \overline{Z}\Bigg\},\\
    g_2&=\dfrac{\partial H_2}{\partial \overline{W}}=\dfrac{1}{2F}\Bigg\{t^2-2-\dfrac{i\gamma_a^2p''(\rho)f(a)}{2a^2}+2t\overline{Z}\sin \theta+\dfrac{\sin \alpha}{2}\overline{Z}- \overline{Z}^2\sin^2\theta\Bigg\},\\
    g_3&=\dfrac{\partial H_2}{\partial \theta}=-\dfrac{\overline{Z}}{2F}\Bigg\{-2t\overline{W}\cos \theta+\dfrac{\cos \alpha}{2}(\overline{Z}- \overline{W})+\overline{W}.\overline{Z}\sin 2\theta \Bigg\},\\
    g_4&=\dfrac{\partial H_2}{\partial t}=\dfrac{H_2}{F}\Bigg\{\bigg[\left(\dfrac{i\gamma_a^2p''(\rho)f(a)}{2a^2}+2(1-t^2)\right)2t-i\gamma_a^2t^2(1-t^2)\left(\gamma_a^2p''(\rho)+\tau'(\rho)\right)\bigg]\\
&~~~~~~~~~~~~~~~~~+3t\overline{Z}+t\overline{W}+\overline{W}.\overline{Z}\sin \theta+i\gamma_a^2t^2\left(\gamma_a^2p''(\rho)+\tau'(\rho)\right)\left(\dfrac{\overline{W}+\overline{Z}}{2}\right)\Bigg\}.
\end{align*}
Therefore, inserting \eqref{eq: 3.73}-\eqref{eq: 3.75} into \eqref{eq: 3.72} and using \eqref{eq: 3.44b}, \eqref{eq: 3.71} one can compute
\begin{align}\label{eq: 3.76}
    \begin{cases}
        \bar{\partial}^+\bar{\partial}^-\overline{W}=\left(\dfrac{5}{2t}+th+tf_1\right)\bar{\partial}^-\overline{W}+G_1,\vspace{0.2 cm}\\
        \bar{\partial}^-\bar{\partial}^+\overline{Z}=\left(\dfrac{5}{2t}+th+tg_1\right)\bar{\partial}^+\overline{Z}+G_2,
    \end{cases}
\end{align}
where 
\begin{align*}
    G_1&=\dfrac{5(\overline{Z}-\overline{W})}{4t^2}+\dfrac{h+f_2}{2}(\overline{Z}-\overline{W})-\dfrac{2H_1+H_2}{2}+\dfrac{t^2f_3\sqrt{1-t^2}(\overline{W}\sin \theta+2t)}{2F}\\
    &~~~~~~~~~+tf_4+(f_2H_2-hH_1)t^2,\\
    G_2&=\dfrac{5(\overline{W}-\overline{Z})}{4t^2}+\dfrac{h+g_2}{2}(\overline{W}-\overline{Z})-\dfrac{H_1+2H_2}{2}+\dfrac{t^2g_3\sqrt{1-t^2}(\overline{Z}\sin \theta-2t)}{2F}\\
    &~~~~~~~~~~~~~~~+tg_4+(g_2H_1-hH_2)t^2.
\end{align*}
We now employ the second order decompositions \eqref{eq: 3.76} to develop the $C^1$ estimates of solutions in the following Lemma:
\begin{l1}\label{l: 4.3}
Let us assume that the conditions \eqref{eq: 2.43} and \eqref{eq: 3.57} are satisfied such that there exists a $C^1$ solution $(W, Z, r, \theta)(t, \psi)$ of the system \eqref{eq: 3.45} with the boundary data \eqref{eq: 2.42} and \eqref{eq: 3.56} in the domain $\Omega_\epsilon$. Then 
\begin{align}
    ||(W, Z, r, \theta)||_{C^1(\Omega_\epsilon)}\leq \dfrac{K_1}{\epsilon^3},
\end{align}
where $K_1$ is a positive constant, independent of $\epsilon$.
\end{l1}
\begin{proof}
From Lemma \ref{l: 4.1}, we know that the functions $W$ and $Z$ are bounded. Therefore, in order to prove this Lemma, it is enough to prove that
\begin{align}
    ||(\overline{W}, \overline{Z}, r, \theta)||_{C^1(\Omega_\epsilon)}\leq \dfrac{K_1}{\epsilon^3},
\end{align}
To prove this, we use the first order characteristic decomposition \eqref{eq: 3.71} and Lemma \ref{l: 4.1} to conclude that
\begin{align}\label{eq: 4.79}
    |\bar{\partial}^+\overline{W}|,~
    |\bar{\partial}^-\overline{Z}|\leq \dfrac{K_1}{t}\leq \dfrac{K_1}{\epsilon} 
\end{align}
for some positive constant $K_1$, independent of $\epsilon$. Also, by Lemma \eqref{l: 4.1} it is easy to observe that the coefficients of $\bar{\partial}^-\overline{W}$ and $\bar{\partial}^+\overline{Z}$ and the functions $G_1$ and $G_2$ in \eqref{eq: 3.76} satisfy
\begin{align*}
    \left|\dfrac{5}{2t}+th+tf_1\right|,~\left|\dfrac{5}{2t}+th+tg_1\right|\leq \dfrac{K_1}{t}, ~|G_1|, ~
    |G_2|\leq \dfrac{K_1}{t^2}.
\end{align*}
Hence, one can integrate the second order decompositions \eqref{eq: 3.76} along the positive and negative characteristic curves to obtain
\begin{align}\label{eq: 4.80}
    |\bar{\partial}^-\overline{W}|,~
    |\bar{\partial}^+\overline{Z}|\leq \dfrac{K_1}{\epsilon}
\end{align}
for some positive constant $K_1$, independent of $\epsilon$. Therefore, noting the expression of scaled normalized derivatives $\bar{\partial}^\pm$, we obtain 
\begin{align*}
    \partial_t=\dfrac{\overline{W}\bar{\partial}^++\overline{Z}\bar{\partial}^-}{\overline{W}+\overline{Z}},~~\partial_\psi=\dfrac{2F\gamma \sqrt{1-t^2}}{ra\gamma_a (\overline{W}+\overline{Z})}.\dfrac{\bar{\partial}^+-\bar{\partial}^-}{t^2},
\end{align*}
which combined with \eqref{eq: 4.79} and \eqref{eq: 4.80} yields the $C^1$ a priori estimates of $(\overline{W}, \overline{Z})$ of the form
\begin{align}
    ||(\overline{W}, \overline{Z})||_{C^1(\Omega_\epsilon)}\leq \dfrac{K_1}{\epsilon^3}.
\end{align}
The estimates of $(r, \theta)$ can be easily derived using the first order decompositions of $r$ and $\theta$ and thus, the proof of the Lemma is completed.
\end{proof}

The existence of a local $C^1$ solution in the neighborhood of the point $T'(\delta, \Tilde{\psi}(\delta))$ can be obtained by using the classical local existence results for boundary value problems to the system of strictly hyperbolic equations \cite{li1985boundary}. Furthermore, utilizing the Lemmas \ref{l: 4.1} and \ref{l: 4.3} we can establish the existence of a global $C^1$ solution for the system \eqref{eq: 3.45} with the boundary data \eqref{eq: 2.42} and \eqref{eq: 3.56} in the domain $\Omega_\epsilon$ by the classical approach of extending  local solution to a larger domain by taking the level sets of $t$ as the Cauchy supports for any $\epsilon>0$. It is noteworthy to see from the system \eqref{eq: 3.45} that the extension step size depends only on the boundary data and the $C^0, C^1$ norms of $(W, Z, r, \theta)$ which are uniformly bounded in $\Omega_\epsilon$. Since the domain $\Omega_\epsilon$ is compact, the extension process can be completed in a finite number of steps. Therefore, by the arbitrariness of $\epsilon$, one can achieve the $C^1$ solution in $\Omega/\{t=0\}$.
\subsection{Regularity of solutions}\label{s: 4.4}
In this subsection, we explore the uniform regularity of solutions up to the degenerate line $t=0$ in the partial hodograph plane. We first derive the uniform boundedness of $L(t, \psi)$ which is uniformly bounded on the boundary curve $\widehat{P'T'}$ by \eqref{eq: 3.40}. Also, noting that $\overline{W}=\dfrac{1}{W}$ and $\overline{Z}=\dfrac{1}{Z}$, we can express $L$ as 
\begin{align}
    L=\dfrac{WZ}{2}.\dfrac{\overline{Z}-\overline{W}}{t}.
\end{align}
Therefore, in order to prove the uniform boundedness of $L$, it is enough to prove that the function $\overline{L}=(\overline{W}-\overline{Z})/t$ is uniformly bounded up to the degenerate line $\widehat{P'D'}$, i.e. $t=0$. A straightforward calculation using \eqref{eq: 3.71} leads us to the equation of $\overline{L}$ of the form
\begin{align}\label{eq: 4.84}
    \begin{cases}
        \bar{\partial}^+\overline{L}=(H_1-H_2)-\dfrac{\bar{\partial}^+\overline{Z}-\bar{\partial}^-\overline{Z}}{t},\\
        \bar{\partial}^-\overline{L}=(H_1-H_2)-\dfrac{\bar{\partial}^+\overline{W}-\bar{\partial}^-\overline{W}}{t}
    \end{cases}
\end{align}
where $H_1$ and $H_2$ are already defined in \eqref{eq: 3.71}.

Therefore, noting the uniform boundedness of $H_1$ and $H_2$ it is necessary to establish the estimates for $(\bar{\partial}^+\overline{W}-\bar{\partial}^-\overline{W})/t$ and $(\bar{\partial}^+\overline{Z}-\bar{\partial}^-\overline{Z})/t$. 

Let us set 
\begin{align}
    U=\bar{\partial}^+\overline{W}-\bar{\partial}^-\overline{W}, ~~V=\bar{\partial}^+\overline{Z}-\bar{\partial}^-\overline{Z}.
\end{align}
Then one can use the commutator relation \eqref{eq: 3.72a} again to obtain the decompositions of $U$ and $V$ of the form
\begin{align}
    \bar{\partial}^+U=\dfrac{\bar{\partial}^-\bar{\lambda}_+-\bar{\partial}^+\bar{\lambda}_-}{\bar{\lambda}_+-\bar{\lambda}_-}U+(\bar{\lambda}_+-\bar{\lambda}_-)(\bar{\partial}^+\overline{W})_\psi
\end{align}
and 
\begin{align}
    \bar{\partial}^-V=\dfrac{\bar{\partial}^-\bar{\lambda}_+-\bar{\partial}^+\bar{\lambda}_-}{\bar{\lambda}_+-\bar{\lambda}_-}V+(\bar{\lambda}_+-\bar{\lambda}_-)(\bar{\partial}^-\overline{Z})_\psi.
\end{align}
Using the same arguments as in the derivation of \eqref{eq: 3.76}, we have
\begin{align}\label{eq: 4.88}
    \begin{cases}
        \bar{\partial}^+U=\left(\dfrac{5}{2t}+t\bar{f}_1\right)U+(2t^2\bar{f}_2-1)\dfrac{V}{2t}+t^2\bar{f}_3,\vspace{0.2 cm}\\
        \bar{\partial}^-V=\left(\dfrac{5}{2t}+t\bar{g}_1\right)V+(2t^2\bar{g}_2-1)\dfrac{U}{2t}+t^2\bar{g}_3,
    \end{cases}
\end{align}
where 
$\bar{f}_1=f_1+h,~~\bar{f}_2=f_2,~~\bar{f}_3=\dfrac{f_3\sqrt{1-t^2}}{2F}[\sin \theta (\overline{Z}-\overline{W})-4t],~~\bar{g}_1=g_1+h,~~\bar{g}_2=g_2,~~\bar{g}_3=\dfrac{g_3\sqrt{1-t^2}}{2F}[\sin \theta (\overline{Z}-\overline{W})-4t]$.

Noting Lemma \ref{l: 4.1} and the uniform boundedness of $f_i$ and $g_i~(i=1, 2, 3)$, it is easy to observe that the functions $\bar{f}_i, \bar{g}_i, (i=1, 2, 3)$ are uniformly bounded in $\Omega$. We further denote
\begin{align}
    \overline{U}=\dfrac{\bar{\partial}^+\overline{W}-\bar{\partial}^- \overline{W}}{t^{\nu}}=\dfrac{U}{t^{\nu}},~~\overline{V}=\dfrac{\bar{\partial}^+\overline{Z}-\bar{\partial}^- \overline{Z}}{t^{\nu}}=\dfrac{V}{t^{\nu}}
\end{align}
for some $\nu\in (0, 2]$.

Then \eqref{eq: 4.88} yields
\begin{align}\label{eq: 3.90}
    \begin{cases}
        \bar{\partial}^+(t^{-\frac{2\nu+1}{2}}\overline{U})=\left(t^2\bar{f}_2-\dfrac{1}{2}\right)t^{-\frac{2\nu+3}{2}}\overline{V}+(t^2\bar{f}_1\overline{U}+t^2\bar{f}_3)t^{-\frac{2\nu+3}{2}},\\
        \bar{\partial}^-(t^{-\frac{2\nu+1}{2}}\overline{V})=\left(t^2\bar{g}_2-\dfrac{1}{2}\right)t^{-\frac{2\nu+3}{2}}\overline{U}+(t^2\bar{g}_1\overline{V}+t^2\bar{g}_3)t^{-\frac{2\nu+3}{2}}.
    \end{cases}
\end{align}

We now use \eqref{eq: 3.90} to prove the uniform boundedness of $\overline{U}$ and $\overline{V}$. Let $A'(0, \psi_1)$ and $A''(0, \psi_2)$ be any two fixed points on the degenerate line  $\widehat{P'D'}$ satisfying $\psi_1>\psi_2$. Then from the points $A'$ and $A''$, we draw positive and negative characteristic curves up to the boundary $\widehat{P'T'}$ at $A_1$ and $A_2$, respectively. Let $\mathcal{D}(\psi_1, \psi_2)\subseteq \Omega$ be the region bounded by the curves $\widehat{A'A''},~\widehat{A'A_1},~~\widehat{A''A_2}$ and $\widehat{A_1A_2}$. Further, let us denote
\begin{align*}
    M'=\underset{\mathcal{D}(\psi_1, \psi_2)}\max \{|\bar{f}_1|, |\bar{g}_1|, |\bar{f}_2|, |\bar{g}_2|, |\bar{f}_3|, |\bar{g}_3|\}.
\end{align*}
\begin{figure}
    \centering
    \includegraphics[width= 4.2 in]{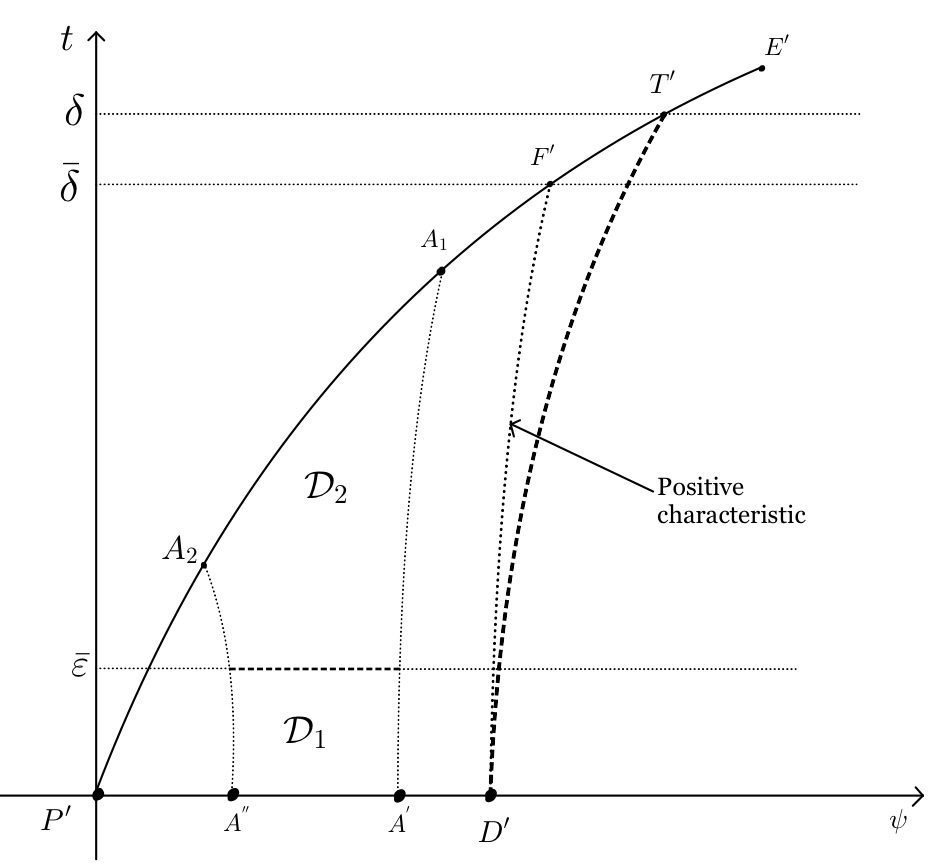}
    \caption{The regions $P'R'D'$, $\mathcal{D}_1$ and $\mathcal{D}_2$}
    \label{fig: 5}
\end{figure}

We then divide the region $\mathcal{D}(\psi_1, \psi_2)$ into two parts, $\mathcal{D}_1:=\mathcal{D}(\psi_1, \psi_2)\cap \{t< \bar{\epsilon}\}$ and $\mathcal{D}_2:=\mathcal{D}(\psi_1, \psi_2)\cap \{t\geq \bar{\epsilon}\}$ (see Figure \ref{fig: 5}), where $\bar{\epsilon}$ is a small number in $(0, \delta]$ satisfying $16\bar{\epsilon}^2M'\leq 1$. Then by the choice of $\bar{\epsilon}$, one easily observe that
\begin{align}\label{eq: 3.91}
    \dfrac{1}{2}+t^2|\bar{f}_2|\leq \dfrac{9}{16},~~t^2|\bar{f}_1|,~t^2|\bar{g}_1|,~t^2|\bar{f}_3|, ~t^2|\bar{g}_3|\leq \dfrac{1}{16}~\forall (t, \psi)\in \overline{\mathcal{D}_1}.
\end{align}
It is easy to observe that we need to check the uniform boundedness of $\overline{U}$ and $\overline{V}$ in the region $\mathcal{D}_1$ only since the degenerate line $t=0$ lies in $\mathcal{D}_1$ and $\overline{U}$ and $\overline{V}$ are uniformly bounded in $\mathcal{D}_2$ and on the boundary $\widehat{A_1A_2}$. Therefore, we now denote
\begin{align}\label{eq: 3.92}
    \widetilde{K}=1+3\max\{\underset{\overline{\mathcal{D}_2}}\max\{|\overline{U}|, |\overline{V}|\}, \underset{\widehat{A_1A_2}}\max\{|\overline{U}|, |\overline{V}|\}\}.
\end{align}
Then we have the following Lemma:
\begin{l1}\label{l: 4.4}
For any point $(t, \psi)\in \mathcal{D}_1$, there holds
\begin{align}
    |\overline{U}|<\widetilde{K}, ~~
    |\overline{V}|<\widetilde{K}.
\end{align}
for any $\nu\in (0, 2]$.
\end{l1}
\begin{proof}
We prove this Lemma by the method of contradiction. Let $B(t_b, \psi_b)$ be any point in $\mathcal{D}_1$ such that point $B$ is the first time that one of $|\overline{U}|$ and $|\overline{V}|$ touches the boundary of $(-\widetilde{K}, \widetilde{K})\times (-\widetilde{K}, \widetilde{K})$. Then from point $B$, we draw a positive and negative characteristic curves up to the upper boundary of $\mathcal{D}_1$ at points $B_1(t_{b_1}, \psi_{b_1})$ and $B_2(t_{b_2}, \psi_{b_2})$, respectively. Without loss of generality, we suppose that $|\overline{U}(t_b, \psi_b)|=\widetilde{K}$ and $|\overline{U}(t, \psi)|\leq \widetilde{K}, ~|\overline{V}(t, \psi)|\leq \widetilde{K}$ hold on the positive characteristic curve $\widehat{BB_1}$. Further, in view of the choice of $\widetilde{K}$, we see
\begin{align*}
    |\overline{U}(t_{b_1}, \psi_{b_1})|, ~|\overline{V}(t_{b_1}, \psi_{b_1})|,~|\overline{U}(t_{b_2}, \psi_{b_2})|,~|\overline{V}(t_{b_2}, \psi_{b_2})|<\dfrac{1}{3}\widetilde{K}.
\end{align*}
Then integrating the equation for $\overline{U}$ in \eqref{eq: 3.90} from $B$ to $B_1$ yields
\begin{align*}
    &t_{b_1}^{-\frac{2\nu+1}{2}}\overline{U}(t_{b_1}, \psi_{b_1})-t_{b}^{-\frac{2\nu+1}{2}}\overline{U}(t_{b}, \psi_{b})=\displaystyle \int_{t_b}^{t_{b_1}}\bigg\{\left(t^2\bar{f}_2-\dfrac{1}{2}\right)t^{-\frac{2\nu+3}{2}}\overline{V}+(t^2\bar{f}_1\overline{U}+t^2\bar{f}_3)t^{-\frac{2\nu+3}{2}}\bigg\}dt,
\end{align*}
which in view of \eqref{eq: 3.91} and \eqref{eq: 3.92} provides
\begin{align*}
    |\overline{U}(t_b, \psi_b)|&\leq t_b^{\frac{2\nu+1}{2}}\bigg\{t_{b_1}^{-\frac{2\nu+1}{2}}|\overline{U}(t_{b_1}, \psi_{b_1})|+\displaystyle \int_{t_b}^{t_{b_1}} \left(\dfrac{9}{16}t^{-\frac{2\nu+3}{2}}|\overline{V}|+\dfrac{1}{16}t^{-\frac{2\nu+3}{2}}|\overline{U}|+\dfrac{1}{16}t^{-\frac{2\nu+3}{2}}\right)dt\bigg\}\\
    &\leq t_b^{\frac{2\nu+1}{2}}\bigg\{\dfrac{1}{3}\widetilde{K}t_{b_1}^{-\frac{2\nu+1}{2}}+\displaystyle \int_{t_b}^{t_{b_1}} \left(\dfrac{9}{16}t^{-\frac{2\nu+3}{2}}\widetilde{K}+\dfrac{1}{16}t^{-\frac{2\nu+3}{2}}\widetilde{K}+\dfrac{1}{16}t^{-\frac{2\nu+3}{2}}\right)dt\bigg\}\\
    &\leq t_b^{\frac{2\nu+1}{2}}\bigg\{\dfrac{1}{3}\widetilde{K}t_{b_1}^{-\frac{2\nu+1}{2}}+\widetilde{K}\displaystyle \int_{t_b}^{t_{b_1}} \dfrac{11}{16}t^{-\frac{2\nu+3}{2}}dt\bigg\}\\
    & = t_b^{\frac{2\nu+1}{2}}\bigg\{\dfrac{1}{3}\widetilde{K}t_{b_1}^{-\frac{2\nu+1}{2}}+\widetilde{K} \dfrac{11}{24}(t_b^{-\frac{2\nu+1}{2}}-t_{b_1}^{-\frac{2\nu+1}{2}})\bigg\}<\dfrac{11}{24}\widetilde{K}<\widetilde{K},
\end{align*}
which contradicts the assumption that $\overline{U}(t_b, \psi_b)=\widetilde{K}$. Therefore, the proof of the Lemma is complete.
\end{proof}

By the arbitrariness of $A'$ and $A''$ on $\widehat{P'D'}$, we can acquire the uniform boundedness of $\overline{U}$ and $\overline{V}$ up to the degenerate line $\widehat{P'D'}$. Then we have the following Lemma:
\begin{l1}\label{l: 4.5}
The functions $\overline{L}$ and $L$ are uniformly bounded up to the degenerate line $\widehat{P'D'}$.
\end{l1}
\begin{proof}
The proof of this Lemma directly follows from Lemma \ref{l: 4.1}, \ref{l: 4.4} and \eqref{eq: 4.84}. Therefore, we omit the details.
\end{proof}
The uniform boundedness of $L$ leads to an important observation that $W=Z$ on the degenerate line $\widehat{P'D'}$. Based on this property, we can develop the uniform regularity of $W, Z$ and $L$ in the following Lemma:
\begin{l1}\label{l: 4.6}
The functions $W, Z$ and $L$ are uniformly $C^{\frac{1}{3}}$ while $r$ and $\theta$ are uniformly $C^{1, \frac{1}{3}}$ continuous in the whole domain $\Omega: P'T'D'$, including the degenerate line $\widehat{P'D'}$.
\end{l1}
\begin{proof}
We first prove the uniform continuity of $W$ at the degenerate line $\widehat{P'D'}$. For any two points $(0, \psi_1)$ and $(0, \psi_2)$ with $\psi_1<\psi_2$, we draw a negative characteristic $l_-$ from $(0, \psi_1)$ and a positive characteristic $l_+$ from $(0, \psi_2)$ and denote the intersection point of $l_-$ and $l_+$ by $(t_m, \psi_m)$ such that the numbers $t_m$ and $\psi_m$ satisfy  
\begin{align}\label{eq: 4.94}
    \psi_m&=\psi_1-\int_{0}^{t_m}\dfrac{ra\gamma_a t^2}{2\gamma FW\sqrt{1-t^2}}dt=\psi_2+\int_{0}^{t_m}\dfrac{ra\gamma_a t^2}{2\gamma FZ\sqrt{1-t^2}}dt.
\end{align}
Now noting the expression of $F(t)$ and the fact that $0<a^2<1$, $0<q<\hat{q}<1$, it is easy to see that there must exist two constants $k_1, k_2>0$ such that $0<\dfrac{k_1}{2}\leq \dfrac{1}{F(t)}\leq \dfrac{k_2}{2}<\infty$. Therefore, using Lemma \ref{l: 4.1}, we must have
\begin{align*}
    \underline{M}t_m^2\leq \dfrac{ra\gamma_a t^2}{2\gamma F\sqrt{1-t^2}}\leq \overline{M} t_m^2.
\end{align*}
for some positive constants $\underline{M}$ and $\overline{M}$.
Then employing \eqref{eq: 4.94} and using Lemma \ref{l: 4.1}, we obtain
\begin{align}\label{eq: 4.95}
    \dfrac{2}{3}\underline{M}t_m^3\leq |\psi_2-\psi_1|\leq \dfrac{2}{3}\overline{M}t_m^3.
\end{align}
Also, according to the uniform boundedness of $L$, we observe by \eqref{eq: 3.45} that there must exists a constant $K>0$ such that $|\partial^+W|\leq K$ and $|\partial^-Z|\leq K$. Therefore, utilizing \eqref{eq: 4.94}, we eventually have
\begin{align*}
|W(0, \psi_2)-W(0, \psi_1)|&=|Z(0, \psi_2)-W(0, \psi_1)|\\
    &\leq |Z(0, \psi_2)-Z(t_m, \psi_m)|+|Z(t_m, \psi_m)-W(t_m, \psi_m)|\\
    &~~~~~~~~+|W(t_m, \psi_m)-W(0, \psi_1)|\\
    &\leq Kt_m+2\underset{\Omega}\max |L|t_m+Kt_m\\
    &\leq (2K+2\underset{\Omega}\max |L|)t_m\leq \widehat{K}|\psi_2-\psi_1|^{\frac{1}{3}}
\end{align*}
for some positive constant $\widehat{K}$, which implies that the function $\overline{W}$ is uniformly $C^{\frac{1}{3}}$ continuous on the degenerate line $\widehat{P'D'}$.

From Lemma \ref{l: 4.3}, it is easy to observe that the function $\overline{W}_t$ is uniformly bounded in the region $\Omega$. Thus if $\psi_1=\psi_2$ then using the boundedness of $\overline{W}_t$, we have
\begin{align*}
    |W(t_1, \psi_1)-W(t_2, \psi_2)|&=
    |W(t_1, \psi_1)-W(t_2, \psi_1)|\\
    &\leq \underset{\Omega}\max|W_t|.|t_2-t_1|\\
    &<K |(t_1, \psi_1)-(t_2, \psi_2)|^\frac{1}{3}.
\end{align*}
For the case $\psi_1<\psi_2$, there are two possible cases:\\\\
\textbf{Case 1.} If $\psi_2-\psi_1\leq t_1$ then we can choose $\nu=\frac{4}{3}$ in Lemma \ref{l: 4.4} and use mean value theorem to obtain
\begin{align*}
   |{W}(t_2, \psi_2)-{W}(t_1, \psi_1)|&\leq |{W}(t_2, \psi_2)-{W}(t_1, \psi_2)|+|{W}(t_1, \psi_2)-{W}(t_1, \psi_1)|\\
    &\leq K|t_2-t_1|+Kt_1^{-\frac{2}{3}}.|\psi_2-\psi_1|\\
    &\leq K|t_2-t_1|+K|\psi_2-\psi_1|^{\frac{1}{3}}\\
    &\leq K|(t_2, \psi_2)-(t_1, \psi_1)|^{\frac{1}{3}} 
\end{align*}
for some uniform constant $K>0$.\\\\
\textbf{Case 2.} For the case $\psi_2-\psi_1>t_1$, we have
\begin{align*}
    |{W}(t_2, \psi_2)-{W}(t_1, \psi_1)|&\leq |{W}(t_2, \psi_2)-{W}(t_1, \psi_2)|+|{W}(t_1, \psi_2)-{W}(0, \psi_2)|\\
    &~~~~~+|{W}(0, \psi_2)-{W}(0, \psi_1)|+|{W}(0, \psi_1)-{W}(t_1, \psi_1)|\\
    &\leq |{W}_t|.|t_2-t_1|+|{W}_t|.t_1+\widehat{K}|\psi_2-\psi_1|^{\frac{1}{3}}+|{W}_t|.t_1\\
    &\leq 2K(|t_2-t_1|+|\psi_2-\psi_1|)+\widehat{K}|\psi_2-\psi_1|^{\frac{1}{3}}\\
    &\leq K'|(t_2, \psi_2)-(t_1, \psi_1)|^{\frac{1}{3}}
\end{align*}
for some uniform constant $K'>0$, which implies that the function ${W}$ is uniformly $C^{\frac{1}{3}}$ continuous in the whole domain $P'T'D'$ including the degenerate line $\widehat{P'D'}$. In a similar manner one can obtain the uniform $C^{\frac{1}{3}}$ continuity of ${Z}$ and $L$ for any two points $(t_1, \psi_1)$ and $(t_2, \psi_2)$ in the domain $P'T'D'$. 

By using the same arguments as above, one can first show that the functions $r(t, \psi)$ and $\theta(t, \psi)$ are uniformly $C^{\frac{1}{3}}$ continuous. Furthermore, we recall from \eqref{eq: 3.44a} that
\begin{align}
    r_t&=-\dfrac{t\sqrt{1-t^2}r \cos \theta}{F(W+Z)},\\
    r_\psi&=\dfrac{\gamma\sqrt{1-t^2}}{a\gamma_a(W+Z)}\left((W+Z)\sin \theta-2\sqrt{1-t^2}\cos \theta L\right),
\end{align}
which means that the function $r(t, \psi)$ is uniformly $C^{1, \frac{1}{3}}$ continuous in the whole domain $P'T'D'$. The same conclusion is also valid for the function $\theta(t, \psi)$ by \eqref{eq: 3.44b}. Hence, the proof of the Lemma is finished.
\end{proof}
Finally, we draw a positive characteristic curve from the point $D'(0, \bar{\psi}(0))$ up to a point $F'(\bar{\delta}, \tilde{\psi}(\bar{\delta}))$ lying on the boundary $\widehat{P'T'}$ (see Figure \ref{fig: 5}) and combine the results of subsections \ref{s: 4.3} and \ref{s: 4.4}to achieve the following theorem:
\begin{t1}\label{t: 4.1}
\textit{Under the assumptions \eqref{eq: 2.43} and \eqref{eq: 3.57}, the system \eqref{eq: 3.45} with boundary data \eqref{eq: 2.42} and \eqref{eq: 3.56} possesses a global smooth solution $\overline{W}(t, r), \overline{Z}(t, r)$ in the entire region $P'F'D'$ bounded by the curves $\widehat{F'D'}$, $\widehat{P'D'}$ and $\widehat{P'F'}$, where $F'$ is the point $(\bar{\delta}, \psi(\bar{\delta}))$ and $\widehat{D'F'}$ is a positive characteristic curve. Further, the solution $(\overline{W}, \overline{Z})(t, \psi)$ and the quantity $\overline{L}(t, \psi)$ are uniformly $C^{\frac{1}{3}}$ continuous up to the degenerate line $\widehat{P'D'}$, i.e. $t=0$.} 
\end{t1}
\section{Solution in the physical plane}
In this section, we recover a global smooth supersonic-sonic solution of the system \eqref{eq: 3.45} in the physical plane using the solutions obtained in Theorem \ref{t: 4.1} in the partial hodograph plane via an inverse transformation.
\subsection{Inversion}
From Theorem \ref{t: 4.1}, we know that the functions $(W, Z, r, \theta)(t, \psi)$ are defined in the whole region $P'F'D'$. Now we proceed to construct the function $x(t, \psi)$ and then prove that the mapping $(t, \psi)\longrightarrow (x, r)$ is a global one-to-one mapping.

Recalling the partial hodograph transformation \eqref{eq: 3.38}, one can easily obtain
\begin{align}
    x_t=\dfrac{t\sqrt{1-t^2}r\sin \theta}{F(W+Z)},~~x_\psi=\dfrac{\sqrt{1-t^2}\gamma}{a\gamma_a (W+Z)}\left((W+Z)\cos \theta+2\sqrt{1-t^2}\sin \theta W\right),
\end{align}
or in other words
\begin{align}\label{eq: 4.99}
    \partial^-x=-\dfrac{tr(t\cos \theta-\sqrt{1-t^2}\sin \theta)}{2FW}.
\end{align}
Then from any point $(\hat{t}, \hat{\psi})$ in the region $P'F'D'$, we draw a negative characteristic curve $\psi=\psi_-(t;\hat{t}, \hat{\psi})~(t\geq \hat{t})$ up to the boundary $\widehat{P'F'}$ at a unique point $(\hat{t}', \Tilde{\psi}(\hat{t}'))$ satisfying
\begin{align}\label{eq: 4.100}
    \begin{cases}
        \dfrac{d\psi_-(t; \hat{t}, \hat{\psi})}{dt}=-\dfrac{r\gamma_a a t^2 }{2\gamma F\sqrt{1-t^2}W}\left(t, \psi_-(t; \hat{t}, \hat{\psi})\right),\\
        \psi_-(\hat{t}; \hat{t}, \hat{\psi})=\hat{\psi}(\hat{t}),~~~\psi_-(\hat{t}'; \hat{t}, \hat{\psi})=\hat{\psi}(\hat{t}')
    \end{cases}
\end{align}
Therefore, we integrate \eqref{eq: 4.99} along the negative characteristic from $\hat{t}'$ to $\hat{t}$ and use \eqref{eq: 4.100} to define the number $x(\hat{t}, \hat{\psi})$ as follows:
\begin{align}
    x(\hat{t}, \hat{\psi})=\hat{x}(\Tilde{\psi}(\hat{t}'))+\displaystyle \int_{\hat{t}}^{\hat{t}'} \dfrac{tr(t\cos \theta-\sqrt{1-t^2}\sin \theta)}{2FW}(t, \psi_-(t; \hat{t}, \hat{\psi}))dt,
\end{align}
where the function $\hat{x}(\psi)$ is defined in \eqref{eq: 3.44d}. Hence, by the arbitrariness of $(\hat{t}, \hat{\psi})$ one can conclude that the function $x=x(t, \psi)$ can be defined in the whole region $P'F'D'$.
\subsubsection{The mapping $(t, \psi)\longrightarrow (x, r)$ is globally injective}
Noting the expressions of $r_t, r_\psi, x_t$ and $x_\psi$, it is straightforward to see that $j:=\dfrac{\partial(x, r)}{\partial(t, \psi)}=\dfrac{t(1-t^2)r\gamma}{a \gamma_a F(t)(W+Z)}\neq 0,$ for $t>0$, which implies that the mapping $(t, \psi)\longrightarrow (x, r)$ is a local one to one mapping. In order to prove that the mapping is globally one to one, including the line $t=0$, we only need to check the strict monotonicity of $\phi$ along the level curves $l^{\epsilon}: 1-\varpi=\epsilon\geq 0$. We prove this by the method of contradiction. Let us assume that there exist two distinct points $(x_1, r_1)$ and $(x_2, r_2)$ in the region $PFD$ such that $t_1=t_2$ and $\psi_1=\psi_2$ which implies that $\cos \omega(x_1, r_1)=\cos \omega(x_2, r_2)$ and $\psi(x_1, r_1)=\psi(x_2, r_2)$ such that both the points $(x_1, r_1)$ and $(x_2, r_2)$ lies on the same level curve $l^\epsilon: 1-\sin \omega=\epsilon\geq 0$. Then we directly compute 
\begin{align*}
    (\phi_x, \phi_r)\cdot (\varpi_r, -\varpi_x)=u\varpi_r-v\varpi_x=\dfrac{\gamma_a(t) \hat{F}_1(t)(W+Z)}{4ar\gamma \sqrt{1-t^2}}>0,
\end{align*}
using the fact that $W, Z>0$ by Lemma \ref{l: 4.1}. Therefore, $\phi$ is monotonically increasing function along each level curve of $l^\epsilon$ which contradicts the assumption that $\phi(x_1, r_1)=\phi(x_2, r_2)$. Hence the mapping is globally injective including the degenerate line $t=0$.
\subsection{Solution of system \eqref{eq: 2.21}}
We now construct the global smooth supersonic solution to system \eqref{eq: 2.21} using the fact that the mapping $(t, \psi)\longrightarrow (x, r)$ is injective so that we can obtain the functions $t=t(x, r)$ and $\psi=\psi(x, r)$ to define the functions
\begin{align}\label{eq: 5.78}
    \theta=\theta(t(x,r), \psi(x, r)),~~\varpi=\sqrt{1-t^2(x, r)}, ~~\forall~(x, r)\in PFD,
\end{align}
where the region $PFD$ is bounded by the curves $\widehat{PF}$, $\widehat{PD}$ and $\widehat{DF}$ such that the curves $\widehat{PD}$ and $\widehat{DF}$ are defined as follows
\begin{align}
\begin{cases}
    \widehat{PD}=\{(x, r)|\varpi(x, r)=1, x\in [x_1, x^*]\},\\
    \widehat{DF}=\{(x, r)|\psi(x,r)=\psi_+(t(x, r); \bar{\delta}, \bar{\psi}(\bar{\delta})), x\in [x^*, x^{**}]\},
    \end{cases}
\end{align}
where $x^*=x(0, \bar{\psi}(0))$ and the number $x^{**}$ satisfies $\Tilde{\psi}(x^{**}, \varphi (x^{**}))=\psi_+(\tilde{t}(x^{**}, \varphi(x^{**})); \bar{\delta}, \tilde{\psi}(\bar{\delta}))$ such that the function $\psi_+$ is the solution of the ODE
\begin{align}
    \begin{cases}
        \dfrac{d\psi_+(t;\bar{\delta}, \Tilde{\psi}(\bar{\delta}))}{dt}=\dfrac{r\gamma_a a t^2}{2\gamma F\sqrt{1-t^2}Z}(t, \psi_+(t; \bar{\delta}, \Tilde{\psi}(\bar{\delta}))), ~~t\in [0, \bar{\delta}]\\
        \psi_+(\bar{\delta};\bar{\delta}, \Tilde{\psi}(\bar{\delta}))=\Tilde{\psi}(\bar{\delta}).
    \end{cases}
\end{align}
We can also get the coordinates of point $D$ and $F$ as $(x^*, r(0, \bar{\psi}(0)))$ and $(x^{**}, \varphi(x^{**}))$, respectively. It is easy to see that the functions $(\theta(x, r), \varpi(x, r))$ defined in \eqref{eq: 5.78} satisfy the boundary condition \eqref{eq: 2.29} by the construction of $(x(t, \psi), r(t, \psi))$. Now we proceed to verify that the function defined in \eqref{eq: 5.78} satisfy the system \eqref{eq: 2.21} in $x-r$ plane.
\subsubsection{Verification of solutions in $x-r$ plane}
By performing a direct calculation, one can yield the following
\begin{align}\label{eq: 5.105}
    \begin{cases}
        \theta_x=\dfrac{t\sin \theta(W-Z)-\sqrt{1-t^2}\cos \theta (W+Z)+\sin^2 \theta}{r}, \vspace{0.2 cm}\\
        \theta_r=\dfrac{t\cos \theta(Z-W)-\sqrt{1-t^2}\sin \theta (W+Z)-\sin \theta \cos \theta}{r}, \vspace{0.2 cm}\\
        \varpi_x=-\dfrac{\hat{F}_1(t)}{4a^2r}\left((W+Z)\sin \theta -2\sqrt{1-t^2} \cos \theta L\right), \vspace{0.2 cm}\\
        \varpi_r=\dfrac{\hat{F}_1(t)}{4a^2r}\left((W+Z)\cos \theta +2\sqrt{1-t^2} \sin \theta L\right).
    \end{cases}
\end{align}
Therefore, using the definition of $\Tilde{\partial}_+$ and \eqref{eq: 5.105}, we obtain
\begin{align*}
    \Tilde{\partial}_+\theta&=r(\cos \alpha\theta_x+\sin \alpha \theta_r)\\
    &= \cos \alpha[\cos \omega \sin \theta (W-Z)-\sin \omega \cos \theta (W+Z)+\sin^2 \theta]\\
    &~~~+\sin \alpha[\cos \omega \cos \theta (W-Z)-\sin \omega \sin \theta (W+Z)-\sin \theta \cos \theta]\\
    &= -\sin \omega \cos \omega(W-Z)-\cos \omega \sin \omega (W+Z)-\varpi \sin \theta\\
    &=-2\varpi \cos \omega W-\varpi \sin \theta.
\end{align*}
and
\begin{align*}
    \Tilde{\partial}_+\varpi&=r(\cos \alpha \varpi_x+\sin \alpha \varpi_r)\\
    &=\dfrac{F_1(\varpi)}{4a^2\gamma}\bigg\{\cos \alpha[-\sin \theta(W+Z)+2\varpi\cos \theta L]+\sin \alpha [\cos \theta(W+Z)+2\varpi \sin \theta L]\bigg\}\\
    &=\dfrac{F_1(\varpi)}{4a^2\gamma}[\varpi(W+Z)+2\varpi \cos \omega L]\\
    &=\dfrac{F_1(\varpi)}{4a^2\gamma}[\varpi(W+Z)+\varpi(W-Z)]\\
    &=2\varpi\dfrac{F_1(\varpi)}{4a^2\gamma}W
\end{align*}
Therefore, we have
\begin{align*}
    \Tilde{\partial}_+\theta+\dfrac{4a^2\gamma\cos \omega}{F_1(\omega)}\Tilde{\partial}_+\varpi&=[-2\varpi\cos \omega W-\varpi \sin \theta]+\dfrac{4a^2\gamma\cos \omega}{F_1(\omega)}\left(2\omega \dfrac{F_1(\omega)}{4a^2\gamma}W\right)\\
    &=-\varpi \sin \theta,
\end{align*}
Hence, the functions $\theta(x, r)$ and $\varpi(x, r)$ satisfy first equation of \eqref{eq: 2.21}. In a similar manner one can prove that $\theta(x, r)$ and $\varpi(x, r)$ satisfy second equation of $\eqref{eq: 2.21}$ and therefore, $\theta(x, r)$ and $\varpi(x, r)$ is the solution of system \eqref{eq: 2.21} with boundary conditions \eqref{eq: 2.29}. 

\subsection{Regularity of angle variables and sonic boundary in the physical plane}
We now discuss the regularity of $\theta(x, r)$ and $\varpi(x, r)$, respectively. Using \eqref{eq: 5.105}, Lemma \ref{l: 4.1} and Lemma \ref{l: 4.5}, it is easy to see that the functions $\theta_x, \theta_r, \varpi_x$ and $\varpi_r$ are uniformly bounded, implying $\theta(x, r)$ and $\varpi(x, r)$ are uniformly Lipschitz continuous. We can actually prove that the functions $\theta_x, \theta_r, \varpi_x$ and $\varpi_r$ are uniformly $C^{\frac{1}{6}}$ continuous. We prove this result in the following Lemma.
\begin{l1}\label{l: 5.1}
Let $f(t, r)$ be a $C^{\frac{1}{3}}$ function defined on the whole region $P'F'D'$. Then if we denote $\bar{f}(x, r)=f(t(x, r), \psi(x, r))$ then the function $\bar{f}(x, y)$ is uniformly $C^{\frac{1}{6}}$ continuous in the whole region $PFD$.
\end{l1}
\begin{proof}
Let $(x', r')$ and $(x''. r'')$ be any two points in $PFD$ and let $(t', \psi')$ and $(t'', \psi'')$ be the images of $(x', r')$ and $(x'', r'')$ in the region $P'F'D'$. Then we evaluate
\begin{align*}
    |\bar{f}(x'', r'')-\bar{f}(x', r')|&=|f(t'', \psi'')-f(t', \psi')|\\
    &\leq K|(t'', \psi'')-(t', \psi')|^{\frac{1}{3}}=K\left(|t''-t'|^2+|\psi''-\psi'|^2\right)^{\frac{1}{6}}
\end{align*}
for some uniform constant $K>0$. Now for the term $|\psi''-\psi'|$, one has
\begin{align*}
    |\psi''-\psi'|=|\phi(x'', r'')-\phi(x', r')|\leq K|(x'', r'')-(x', r')|,
\end{align*}
by the uniform Lipschitz continuity of $\phi(x, r)$. 

Again, for the term $|t''-t'|^2$, we have
\begin{align*}
    |t''-t'|^2&\leq |t''-t'|.|t''+t'|=|t''^2-t'^2|\\
    &\leq |(1-\varpi^2(x'', r''))-(1-\varpi^2(x', r'))|\\
    &\leq |\varpi^2(x'', r'')-\varpi^2(x', r')|\leq 2|\varpi(x'', r'')-\varpi(x', r')|\\
    &\leq K|(x'', r'')-(x', r')|
\end{align*}
in view of the uniform Lipschitz continuity of $\varpi$. Hence, we directly compute
\begin{align}
    |\bar{f}(x'', r'')-\bar{f}(x', r')|\leq K|(x'', r'')-(x', r')|^{\frac{1}{6}},
\end{align}
which implies that the function $\bar{f}$ is uniformly $C^{\frac{1}{6}}$ continuous. Therefore, Lemma is proved.
\end{proof}
Hence, by using the Lemma \ref{l: 4.5} and Lemma \ref{l: 5.1}, it is easy to see that the functions $(W, Z, L, \theta)(t(x, r), \psi(x, r))$ are uniformly $C^{\frac{1}{6}}$ continuous in the whole region $PFD$ up to the sonic boundary $\widehat{PD}$. Therefore, using \eqref{eq: 5.105} we observe that the functions $\varpi(x, r)$ and $\theta(x, r)$ are uniformly $C^{1, \frac{1}{6}}$ continuous in the domain $PFD$ up to the sonic curve $\widehat{PD}$. 

Moreover, we directly apply \eqref{eq: 5.105}, Lemma \ref{l: 4.1} and Lemma \ref{l: 4.5} to observe that
\begin{align*}
    0<\widetilde{k}\leq (\varpi_x)^2+(\varpi_r)^2=\dfrac{(F_1(\varpi))^2}{16a^4r^2}\left((W+Z)^2+4\varpi^2 L^2\right)\leq \widetilde{K}<\infty,
\end{align*}
for some positive constants $\widetilde{K}$ and $\widetilde{K}$. This fact implies that the level curve $\varpi(x, r)=\epsilon\geq 0$ is $C^1$ continuous. Furthermore, due to Lemmas \ref{l: 4.5} and \ref{l: 5.1}, the level curve $\varpi(x, r)$ and eventually the sonic curve $\widehat{PD}$ are actually $C^{1, \frac{1}{6}}$ continuous. 
\subsection{$\widehat{DF}$ is a negative characteristic curve}
Since $\widehat{D'F'}$ is a positive characteristic curve so in order to prove that $\widehat{DF}$ is a negative characteristic curve, it is enough to prove that the mapping $(t, \psi)\longrightarrow (x, r)$ transforms a positive characteristic curve in $t-\psi$ plane into a negative characteristic curve in $x-r$ plane. We prove this in the following Lemma.
\begin{l1}\label{l: 5.2}
A positive characteristic curve in $(t, \psi)$ plane is transformed into a negative characteristic curve in $(x, r)$ plane under the transformation \eqref{eq: 5.78}.
\end{l1}
\begin{proof}
To prove it, we differentiate the equality $\psi(x, r)=\psi_+(t(x, r))$ with respect to $x$ and use the fact $\psi_+'(t)=\lambda$ to get
\begin{align}\label{eq: 5.86}
    \dfrac{dr}{dx}=-\dfrac{\psi_x-\lambda t_x}{\psi_r-\lambda t_r}=-\dfrac{\lambda \sqrt{1-t^2}\omega_x-\phi_x}{\lambda \sqrt{1-t^2}\omega_r-\phi_r}.
\end{align}
Then by exploiting \eqref{eq: 5.105} in \eqref{eq: 5.86} and applying \eqref{eq: 5.78} yields
\begin{align*}
    \dfrac{dr}{dx}&=\dfrac{-2\sqrt{1-t^2}Z\cos \theta+t(W+Z)\sin \theta-2\sqrt{1-t^2}\cos \theta Lt}{2\sqrt{1-t^2}Z\sin \theta+t(W+Z)\cos \theta+2\sqrt{1-t^2}\sin \theta Lt}\\
    &=\dfrac{\cos \omega \sin \theta-\sin \omega \cos \theta}{\cos \omega \cos \theta+\sin \omega \sin \theta}=\dfrac{\sin \beta}{\cos \beta}=\Lambda_-,
\end{align*}
which implies that the curves defined by the equality $\psi(x, r)=\psi_+(t(x, r))$ are negative characteristics in the $x-r$ plane. 
\end{proof}
Now since $\dfrac{da}{dq}<0$ and $\dfrac{d\gamma}{dq}>0$, then we must have $\dfrac{da}{d\gamma}<0$ or in other words $\gamma:=\gamma(a)$. Then we combine $a=a(\varpi)$ and \eqref{eq: 5.78} together with $u=\dfrac{a\gamma_a \cos \theta}{\gamma \sin \omega}, v=\dfrac{a\gamma_a \sin \theta}{\gamma \sin \omega}$ to define the functions $(a, u, v)(x, r)$ such that
\begin{align}
    a(x, r)=a(\varpi(x, r)),~u=\dfrac{a(\varpi(x, r))\gamma_{a(\varpi(x, r))} \cos \theta(x, r)}{\gamma(a(\varpi(x, r))) \varpi(x, r)},~v=\dfrac{a(\varpi(x, r))\gamma_{a(\varpi(x, r))} \sin \theta(x, r)}{\gamma(a(\varpi(x, r))) \varpi(x, r)}
\end{align}
is the classical solution of \eqref{eq: 2.9}.

In view of assumptions \eqref{eq: 2.31a}, \eqref{eq: 2.31}, Theorem \ref{t: 4.1}, Lemma \ref{l: 5.1} and Lemma \ref{l: 5.2}, we have the following result:
\begin{t1}
\textit{Let $\widehat{PE}: r=\varphi(x)$ is an increasing and concave smooth streamline of a 3-D steady axisymmetric isentropic irrotational relativistic flow such that the Mach number $M$ increases along $\widehat{PE}$ with $M=1$ at the point $P$ and $\varphi'$ and $\varpi$ satisfy \eqref{eq: 2.31a} and \eqref{eq: 2.31}. Then there exists a smooth sonic curve $\widehat{PD}$ and a negative characteristic curve $\widehat{DF}$ such that the boundary value problem \eqref{eq: 2.21}-\eqref{eq: 2.29} has a smooth supersonic solution $(\theta, \varpi)$ in the region $PFD$, where $F$ is a point lying on the streamline $\widehat{PE}$. Furthermore, solution $(\theta, \varpi)(x, r)$ is uniformly $C^{1, \frac{1}{6}}$ continuous in the whole region $PFD$ while the sonic curve $\widehat{PD}$ is $C^{1, \frac{1}{6}}$ continuous.} 
\end{t1}
\section{Conclusions}\label{9} 
In this article, we considered three-dimensional axisymmetric steady isentropic relativistic Euler equations with a general convex pressure and proved the global existence and regularity of solution of a supersonic-sonic patch arising in the modified Frankl problem. Using the characteristic decompositions of angle variables and a partial hodograph transformation, we were able to prove that solution is uniformly $C^{1, \frac{1}{6}}$ continuous. Moreover, we proved that the sonic boundary is $C^{1, \frac{1}{6}}$ continuous. The study of such supersonic-sonic patch problems is quite crucial in the context of transonic flows. Here we constructed solution up to a negative characteristic curve $\widehat{DF}$. However, in the future, we will try to construct a global smooth supersonic solution of the modified Frankl problem for relativistic Euler equations with arbitrary equation of state up to the positive characteristic curve $\widehat{EO}$ by solving a free boundary value problem and using the symmetry of the airfoil.
\section*{Acknowledgments}
\textit{The first author (RB) gratefully acknowledges the research support from the University Grant Commission, Government of India.  The second author (TRS) would like to
thank SERB, DST, India (Ref. No. MTR/2019/001210) for its financial support through the MATRICS grant.}
\biboptions{sort&compress}
\bibliographystyle{elsarticle-num}
\bibliography{Reference}%

\begin{thebibliography}{10}
\expandafter\ifx\csname url\endcsname\relax
  \def\url#1{\texttt{#1}}\fi
\expandafter\ifx\csname urlprefix\endcsname\relax\def\urlprefix{URL }\fi
\expandafter\ifx\csname href\endcsname\relax
  \def\href#1#2{#2} \def\path#1{#1}\fi

\bibitem{courant0}
R.~Courant, K.~O. Friedrichs, Supersonic flow and shock waves, Vol.~21,
  Springer Science \& Business Media, 1999.

\bibitem{bers2016mathematical}
L.~Bers, Mathematical aspects of subsonic and transonic gas dynamics, Courier
  Dover Publications, 2016.

\bibitem{kuz2003boundary}
A.~G. Kuz'min, Boundary value problems for transonic flow, John Wiley \& Sons,
  2003.

\bibitem{shapiro1953dynamics}
A.~H. Shapiro, The dynamics and thermodynamics of compressible fluid flow, New
  York: Ronald Press (1953).

\bibitem{li1998two}
J.~Li, T.~Zhang, S.~Yang, The two-dimensional \relax{R}iemann problem in gas
  dynamics, Vol.~98, CRC Press, 1998.

\bibitem{li2009interaction}
J.~Li, Y.~Zheng, Interaction of rarefaction waves of the two-dimensional
  self-similar \relax{E}uler equations, Archive for rational mechanics and
  analysis 193~(3) (2009) 623--657.

\bibitem{li2011characteristic}
J.~Li, Z.~Yang, Y.~Zheng, Characteristic decompositions and interactions of
  rarefaction waves of 2-\relax{D} \relax{E}uler equations, Journal of
  Differential Equations 250~(2) (2011) 782--798.

\bibitem{gilbarg1955uniqueness}
D.~Gilbarg, J.~Serrin, Uniqueness of axially symmetric subsonic flow past a
  finite body, Journal of Rational Mechanics and Analysis 4 (1955) 169--175.

\bibitem{xie2007global}
C.~Xie, Z.~Xin, Global subsonic and subsonic-sonic flows through infinitely
  long nozzles, Indiana University mathematics journal (2007) 2991--3023.

\bibitem{xie2010global}
C.~Xie, Z.~Xin, Global subsonic and subsonic-sonic flows through infinitely
  long axially symmetric nozzles, Journal of Differential Equations 248~(11)
  (2010) 2657--2683.

\bibitem{chen2016subsonic}
G.-Q. Chen, F.-M. Huang, T.-Y. Wang, Subsonic-sonic limit of approximate
  solutions to multidimensional steady \relax{E}uler equations, Archive for
  Rational Mechanics and Analysis 219~(2) (2016) 719--740.

\bibitem{wang2019smooth}
C.~Wang, Z.~Xin, Smooth transonic flows of \relax{Meyer type in de Laval
  nozzles}, Archive for Rational Mechanics and Analysis 232~(3) (2019)
  1597--1647.

\bibitem{wang2021regular}
C.~Wang, Z.~Xin, Regular subsonic-sonic flows in general nozzles, Advances in
  Mathematics 380 (2021) 107578.

\bibitem{zhang2014sonic}
T.~Zhang, Y.~Zheng, Sonic-supersonic solutions for the steady \relax{E}uler
  equations, Indiana University Mathematics Journal (2014) 1785--1817.

\bibitem{hu2019sonic}
Y.~Hu, T.~Li, Sonic-supersonic solutions for the two-dimensional pseudo-steady
  full \relax{E}uler equations, Kinetic \& Related Models 12~(6) (2019) 1197.

\bibitem{hu2020sonic}
Y.~Hu, J.~Li, Sonic-supersonic solutions for the two-dimensional steady full
  \relax{E}uler equations, Archive for Rational Mechanics and Analysis 235~(3)
  (2020) 1819--1871.

\bibitem{li2019degenerate}
F.~Li, Y.~Hu, \relax{On a degenerate mixed-type boundary value problem to the
  2-D steady Euler equations}, Journal of Differential Equations 267~(11)
  (2019) 6265--6289.

\bibitem{hu2020global}
Y.~Hu, J.~Li, On a global supersonic-sonic patch characterized by \relax{2-D}
  steady full \relax{E}uler equations, Advances in Differential Equations
  25~(5/6) (2020) 213--254.

\bibitem{du2011subsonic}
L.~Du, Z.~Xin, W.~Yan, Subsonic flows in a multi-dimensional nozzle, Archive
  for rational mechanics and analysis 201~(3) (2011) 965--1012.

\bibitem{du2014steady}
L.~Du, C.~Xie, Z.~Xin, Steady subsonic ideal flows through an infinitely long
  nozzle with large vorticity, Communications in Mathematical Physics 328~(1)
  (2014) 327--354.

\bibitem{chen2016two}
C.~Chen, L.~Du, C.~Xie, Z.~Xin, Two dimensional subsonic \relax{E}uler flows
  past a wall or a symmetric body, Archive for Rational Mechanics and Analysis
  221~(2) (2016) 559--602.

\bibitem{wang2013degenerate}
C.~Wang, Z.~Xin, On a degenerate free boundary problem and continuous
  subsonic--sonic flows in a convergent nozzle, Archive for Rational Mechanics
  and Analysis 208~(3) (2013) 911--975.

\bibitem{hu2021sonic}
Y.~Hu, J.~Chen, Sonic-supersonic solutions to a mixed-type boundary value
  problem for the two-dimensional full \relax{E}uler equations, SIAM Journal on
  Mathematical Analysis 53~(2) (2021) 1579--1629.

\bibitem{chen2007two}
G.-Q. Chen, C.~M. Dafermos, M.~Slemrod, D.~Wang, On two-dimensional
  sonic-subsonic flow, Communications in mathematical physics 271~(3) (2007)
  635--647.

\bibitem{morawetz1964non}
C.~S. Morawetz, Non-existence of transonic flow past a profile, Communications
  on Pure and Applied Mathematics 17~(3) (1964) 357--367.

\bibitem{frankl1950formation}
F.~Frankl, On the formation of shock waves in subsonic flows with local
  supersonic velocities, Prikladnaya Matematika I Mekhanika 11~(NACA-TM-1251)
  (1950).

\bibitem{morawetz1954uniqueness}
C.~S. Morawetz, A uniqueness theorem for \relax{F}rankl's problem,
  Communications on Pure and Applied Mathematics 7~(4) (1954) 697--703.

\bibitem{cook1978uniqueness}
L.~P. Cook, A uniqueness proof for a transonic flow problem, Indiana University
  Mathematics Journal 27~(1) (1978) 51--71.

\bibitem{lighthill1945new}
M.~J. Lighthill, A new method of two-dimensional aerodynamic design (1945).

\bibitem{hassan1981transonic}
A.~Hassan, H.~Sobieczky, Transonic airfoils with a given pressure distribution,
  in: 14th Fluid and Plasma Dynamics Conference, p. 1235.

\bibitem{henne1981inverse}
P.~Henne, Inverse transonic wing design method, Journal of Aircraft 18~(2)
  (1981) 121--127.

\bibitem{volpe1981role}
G.~Volpe, R.~MELNICK, The role of constraints in the inverse design problem for
  transonic airfoils, in: 14th Fluid and Plasma Dynamics Conference, 1981, p.
  1233.

\bibitem{volpe1986design}
G.~Volpe, R.~Melnik, The design of transonic aerofoils by a well-posed inverse
  method, International journal for numerical methods in engineering 22~(2)
  (1986) 341--361.

\bibitem{stanitz1988review}
J.~D. Stanitz, A review of certain inverse methods for the design of ducts with
  2-or 3-dimensional potential flow (1988).

\bibitem{labrujere1993computational}
T.~E. Labrujere, J.~Slooff, Computational methods for the aerodynamic design of
  aircraft components, Annual Review of Fluid Mechanics 25~(1) (1993) 183--214.

\bibitem{obayashi1996genetic}
S.~Obayashi, S.~Takanashi, Genetic optimization of target pressure
  distributions for inverse design methods, AIAA journal 34~(5) (1996)
  881--886.

\bibitem{kuz2001solvability}
A.~Kuz'min, Solvability of a problem for transonic flow with a local supersonic
  region, Nonlinear Differential Equations and Applications NoDEA 8~(3) (2001)
  299--321.

\bibitem{kuz2004modified}
A.~G. Kuz’min, A modified \relax{Frankl-Morawetz problem on a transonic flow
  past an airfoil}, Differential Equations 40~(10) (2004) 1455--1460.

\bibitem{husonic2021}
Y.~Hu, J.~Li, On a supersonic-sonic patch arising from the \relax{F}rankl
  problem in transonic flows, Communications on Pure \& Applied Analysis 20~(7
  \& 8) (2021) 2643--2663.

\bibitem{hu2022supersonic}
Y.~Hu, On a supersonic-sonic patch in the three-dimensional steady axisymmetric
  transonic flows, SIAM Journal on Mathematical Analysis 54~(2) (2022)
  1515--1542.

\bibitem{luan2018two}
L.~Luan, J.~Chen, J.~Liu, \relax{Two dimensional relativistic Euler equations
  in a convex duct}, Journal of Mathematical Analysis and Applications 461~(2)
  (2018) 1084--1099.

\bibitem{li2005global}
Y.~Li, D.~Feng, Z.~Wang, \relax{Global entropy solutions to the relativistic
  Euler equations for a class of large initial data}, Zeitschrift f{\"u}r
  angewandte Mathematik und Physik ZAMP 56~(2) (2005) 239--253.

\bibitem{chen2018boundary}
J.~Chen, G.~Lai, J.~Zhang, \relax{Boundary value problems for the 2D steady
  relativistic Euler equations with general equation of state}, Nonlinear
  Analysis 175 (2018) 56--72.

\bibitem{chen2004stability}
G.-Q. Chen, Y.~Li, \relax{Stability of Riemann solutions with large oscillation
  for the relativistic Euler equations}, Journal of Differential Equations
  202~(2) (2004) 332--353.

\bibitem{li2011semi}
M.~Li, Y.~Zheng, Semi-hyperbolic patches of solutions to the two-dimensional
  \relax{E}uler equations, Archive for Rational Mechanics and Analysis 201~(3)
  (2011) 1069--1096.

\bibitem{barthwal2022existence}
R.~Barthwal, T.~Raja~Sekhar, On the existence and regularity of solutions of
  semihyperbolic patches to 2-\relax{D Euler equations with van der Waals gas},
  Studies in Applied Mathematics 148~(2) (2022) 543--576.

\bibitem{fan2022sonic}
Y.~Fan, L.~Guo, Y.~Hu, S.~You, Sonic-supersonic solutions to a degenerate
  \relax{Cauchy--Goursat problem for 2D relativistic Euler equations},
  Zeitschrift f{\"u}r angewandte Mathematik und Physik 73~(1) (2022) 1--24.

\bibitem{lai2015centered}
G.~Lai, W.~Sheng, Centered wave bubbles with sonic boundary of pseudosteady
  \relax{G}uderley \relax{M}ach reflection configurations in gas dynamics,
  Journal de Math{\'e}matiques Pures et Appliqu{\'e}es 104~(1) (2015) 179--206.

\bibitem{sheng2018interaction}
W.~Sheng, S.~You, Interaction of a centered simple wave and a planar
  rarefaction wave of the two-dimensional \relax{E}uler equations for
  pseudo-steady compressible flow, Journal de Mathematiques Pures et Appliquees
  114~(9) (2018) 29--50.

\bibitem{lai2014characteristic}
G.~Lai, C.~Shen, \relax{Characteristic decompositions and boundary value
  problems for two-dimensional steady relativistic Euler equations},
  Mathematical Methods in the Applied Sciences 37~(1) (2014) 136--147.

\bibitem{hu2021degenerate}
Y.~Hu, F.~Li, On a degenerate hyperbolic problem for the 3-\relax{D} steady
  full \relax{E}uler equations with axial-symmetry, Advances in Nonlinear
  Analysis 10~(1) (2021) 584--615.

\bibitem{li1985boundary}
D.~Li, W.~Yu, Boundary value problems for quasilinear hyperbolic systems, Duke
  University, 1985.

\end{thebibliography}


\begin{thebibliography}{10}

\bibitem{Hirt1974}
Hirt CW, Amsden AA, Cook JL. An arbitrary {L}agrangian-{E}ulerian computing
  method for all flow speeds.  {\it J {C}omput {P}hys. }1974;14(3):227--253.

\bibitem{Liska2010}
Liska R, Shashkov M, Vachal P, Wendroff B. Optimization-based synchronized
  flux-corrected conservative interpolation (remapping) of mass and momentum
  for arbitrary {L}agrangian-{E}ulerian methods.  {\it J {C}omput {P}hys.
  }2010;229(5):1467--1497.

\bibitem{Taylor1937}
Taylor GI, Green AE. Mechanism of the production of small eddies from large
  ones.  {\it P {R}oy {S}oc {L}ond {A} {M}at. }1937;158(895):499--521.
\newblock \url{https://doi.org/10.1098/rspa.1937.0036},
  \url{http://rspa.royalsocietypublishing.org/content/158/895/499}.

\bibitem{Knupp1999}
Knupp PM. Winslow smoothing on two-dimensional unstructured meshes.  {\it Eng
  {C}omput. }1999;15:263--268.

\bibitem{Kamm2000}
Kamm J. {\it Evaluation of the {S}edov-von {N}eumann-{T}aylor blast wave
  solution. } Technical {R}eport LA-UR-00-6055: Los {A}lamos {N}ational
  {L}aboratory; 2000.

\bibitem{Kucharik2003}
Kucharik M, Shashkov M, Wendroff B. An efficient linearity-and-bound-preserving
  remapping method.  {\it J {C}omput {P}hys. }2003;188(2):462--471.

\bibitem{Blanchard2015}
Blanchard G, Loubere R. {\it High-Order {C}onservative {R}emapping with a
  posteriori {MOOD} stabilization on polygonal meshes. }
  \url{https://hal.archives-ouvertes.fr/hal-01207156}, the {HAL} {O}pen
  {A}rchive, hal-01207156. Accessed January 13, 2016; 2015.

\bibitem{Burton2013}
Burton DE, Kenamond MA, Morgan NR, Carney TC, Shashkov MJ. An intersection
  based {ALE} scheme {(xALE)} for cell centered hydrodynamics {(CCH)}.  In:
  Talk at {M}ultimat 2013, {I}nternational {C}onference on {N}umerical
  {M}ethods for {M}ulti-{M}aterial {F}luid {F}lows; September 2--6, 2013; San
  {F}rancisco.
\newblock LA-UR-13-26756.2.

\bibitem{Berndt2011}
Berndt M, Breil J, Galera S, Kucharik M, Maire PH, Shashkov M. Two-step hybrid
  conservative remapping for multimaterial arbitrary {L}agrangian-{E}ulerian
  methods.  {\it J {C}omput {P}hys. }2011;230(17):6664--6687.

\bibitem{Kucharik2012}
Kucharik M, Shashkov M. One-step hybrid remapping algorithm for multi-material
  arbitrary {L}agrangian-{E}ulerian methods.  {\it J {C}omput {P}hys.
  }2012;231(7):2851--2864.

\bibitem{Breil2015}
Breil J, Alcin H, Maire PH. A swept intersection-based remapping method for
  axisymmetric {ReALE} computation.  {\it Int {J} {N}umer {M}eth {F}l.
  }2015;77(11):694--706.
\newblock Fld.3996.

\bibitem{Barth1997}
Barth TJ. Numerical methods for gasdynamic systems on unstructured meshes.  In:
   Kroner D, Rohde C, Ohlberger M, eds. {\it An {I}ntroduction to {R}ecent
  {D}evelopments in {T}heory and {N}umerics for {C}onservation {L}aws,
  {P}roceedings of the {I}nternational {S}chool on {T}heory and {N}umerics for
  {C}onservation {L}aws}, Lecture {N}otes in {C}omputational {S}cience and
  {E}ngineering. Berlin: Springer 1997.
\newblock ISBN 3-540-65081-4.

\bibitem{Lauritzen2011}
Lauritzen P, Erath C, Mittal R. On simplifying `incremental remap'-based
  transport schemes.  {\it J {C}omput {P}hys. }2011;230(22):7957--7963.

\bibitem{Klima2017}
Klima M, Kucharik M, Shashkov M. Local error analysis and comparison of the
  swept- and intersection-based remapping methods.  {\it Commun {C}omput
  {P}hys. }2017;21(2):526--558.

\bibitem{Dukowicz2000}
Dukowicz JK, Baumgardner JR. Incremental remapping as a transport/advection
  algorithm.  {\it J {C}omput {P}hys. }2000;160(1):318--335.

\bibitem{Kucharik2011}
Kucharik M, Shashkov M. Flux-based approach for conservative remap of
  multi-material quantities in {2D} arbitrary {L}agrangian-{E}ulerian
  simulations.  In:  Fo\v{r}t J, F{\"{u}}rst J, Halama J, Herbin R, Hubert F,
  eds. {\it Finite {V}olumes for {C}omplex {A}pplications {VI} {P}roblems \&
  {P}erspectives},  Springer {P}roceedings in {M}athematics, vol. 1: Springer
  2011 (pp. 623--631).

\bibitem{Kucharik2014}
Kucharik M, Shashkov M. Conservative multi-material remap for staggered
  multi-material arbitrary {L}agrangian-{E}ulerian methods.  {\it J {C}omput
  {P}hys. }2014;258:268--304.

\bibitem{Loubere2005}
Loubere R, Shashkov M. A subcell remapping method on staggered polygonal grids
  for arbitrary-{L}agrangian-{E}ulerian methods.  {\it J {C}omput {P}hys.
  }2005;209(1):105--138.

\bibitem{Caramana1998}
Caramana EJ, Shashkov MJ. Elimination of artificial grid distortion and
  hourglass-type motions by means of {L}agrangian subzonal masses and
  pressures.  {\it J {C}omput {P}hys. }1998;142(2):521--561.

\bibitem{Hoch2009}
Hoch P. {\it An arbitrary {L}agrangian-{E}ulerian strategy to solve
  compressible fluid flows. } Technical {R}eport: CEA; 2009.
\newblock HAL: hal-00366858.
  https://hal.archives-ouvertes.fr/docs/00/36/68/58/PDF/ale2d.pdf. Accessed
  January 13, 2016.

\bibitem{Shashkov1996}
Shashkov M. {\it Conservative {F}inite-{D}ifference {M}ethods on {G}eneral
  {G}rids}.
\newblock Boca Raton, Florida: CRC {P}ress; 1996.
\newblock ISBN 0-8493-7375-1.

\bibitem{Benson1992}
Benson DJ. Computational methods in {L}agrangian and {E}ulerian hydrocodes.
  {\it Comput {M}ethod {A}ppl {M}. }1992;99(2--3):235--394.

\bibitem{Margolin2003}
Margolin LG, Shashkov M. Second-order sign-preserving conservative
  interpolation (remapping) on general grids.  {\it J {C}omput {P}hys.
  }2003;184(1):266--298.

\bibitem{Kenamond2013}
Kenamond MA, Burton DE. Exact intersection remapping of multi-material
  domain-decomposed polygonal meshes.  In: Talk at {M}ultimat 2013,
  {I}nternational {C}onference on {N}umerical {M}ethods for {M}ulti-{M}aterial
  {F}luid {F}lows; September 2--6, 2013; San {F}rancisco.
\newblock LA-UR-13-26794.

\bibitem{Dukowicz1984}
Dukowicz J. Conservative rezoning (remapping) for general quadrilateral meshes.
   {\it J {C}omput {P}hys. }1984;54(3):411--424.

\bibitem{Margolin2002}
Margolin LG, Shashkov M. {\it Second-order sign-preserving remapping on general
  grids. } Technical Report LA-UR-02-525: Los {A}lamos {N}ational {L}aboratory;
  2002.

\bibitem{Mavriplis2003}
Mavriplis DJ. Revisiting the least-squares procedure for gradient
  reconstruction on unstructured meshes.  In: AIAA 2003-3986. 16th {AIAA}
  {C}omputational {F}luid {D}ynamics {C}onference; June 23--26, 2003; Orlando,
  {F}lorida.

\bibitem{Scovazzi2008}
Scovazzi G, Love E, Shashkov M. Multi-scale {L}agrangian shock hydrodynamics on
  {Q1/P0} finite elements: {T}heoretical framework and two-dimensional
  computations.  {\it Comput {M}ethod {A}ppl {M}. }2008;197(9--12):1056--1079.

\end{thebibliography}
\end{document}